\documentclass[a4paper,11pt,leqno]{article}
\usepackage{amsfonts,amssymb,amsmath,amsthm}
\usepackage{fullpage}
\usepackage{epic,eepic}
\usepackage{mathrsfs}%adds \mathscr{} fons

\newtheorem{theorem}{Theorem}[section]
\newtheorem{proposition}[theorem]{Proposition}
\newtheorem{corollary}[theorem]{Corollary}
\newtheorem{lemma}[theorem]{Lemma}
\theoremstyle{remark}
\newtheorem{remark}[theorem]{Remark}
\numberwithin{equation}{section}

\newcommand{\N}{{\mathbb N}}
\newcommand{\R}{{\mathbb R}}
\newcommand{\C}{{\mathcal C}}
\newcommand{\E}{{\mathcal E}}
\newcommand{\D}{{\mathcal D}}
\newcommand{\eps}{{\varepsilon}}

\begin{document}

\title{Positive solutions to singular semilinear elliptic equations \\
with critical potential on cone--like domains}
\author{Vitali Liskevich, Sofya Lyakhova and Vitaly Moroz\medskip\\
School of Mathematics\\
University of Bristol\\
Bristol BS8 1TW\\
United Kingdom\smallskip\\
{\small E-mail: {\tt V.Liskevich$\,|$S.Lyakhova$\,|$V.Moroz@bristol.ac.uk}}}

\date{}

\maketitle

\begin{abstract}
We study the existence and nonexistence of positive (super-)\,solutions
to a singular semilinear elliptic equation
$$-\nabla\cdot(|x|^A\nabla u)-B|x|^{A-2}u=C|x|^{A-\sigma}u^p$$
in cone--like domains of $\R^N$ ($N\ge 2$), for the full range of
parameters $A,B,\sigma,p\in\R$ and $C>0$. We provide a complete
characterization of the set of $(p,\sigma)\in\R^2$ such that the
equation has no positive (super-)\,solutions, depending on the
values of $A,B$ and the principle Dirichlet eigenvalue of the
cross--section of the cone.

The proofs are based on the explicit construction of appropriate barriers and involve
the analysis of asymptotic behavior of super-harmonic functions associated
to the Laplace operator with critical potentials, Phragmen--Lindel\"of type comparison arguments
and an improved version of Hardy's inequality in cone--like domains.
\end{abstract}

\begin{small}
\tableofcontents
\end{small}

%\newpage

\section{Introduction and main results}

We study the existence and nonexistence of positive (super)\,solutions
to the singular semilinear elliptic equation with critical potential
\begin{equation}\label{e:MAIN}
-\nabla\cdot(|x|^A\nabla u)-B|x|^{A-2}u=C|x|^{A-\sigma}u^p\quad\mbox{in }\C_\Omega^\rho.
\end{equation}
Here $A,B\in\R$, $C>0$ and $(p,\sigma)\in\R^2$. By
$\C^\rho_\Omega\subset\R^N$ ($N\ge 2$) we denote the cone-like domain defined by
$$\C^\rho_\Omega=\{(r,\omega)\in \R^N:\omega\in\Omega,\,r>\rho\},$$
where $\rho>0$, $(r,\omega)$ are the polar coordinates in $\R^N$
and $\Omega\subseteq S^{N-1}$ is a subdomain (connected open
subset) of the unit sphere $S^{N-1}$ in $\R^N$. Note that we do
not prescribe any boundary conditions in \eqref{e:MAIN}.
A nonegative \textit{super-solution} to \eqref{e:MAIN} in a domain
$G\subseteq\R^N$ is an $0\leq u\in H^1_{loc}(G)$ such that
\begin{equation}\label{e:DefSol}
\int_{G}\nabla u\cdot\nabla\varphi\,|x|^A dx- B\int_{G}
u\varphi\,|x|^{A-2}dx\,
\geq\,C\int_{G}u^p\varphi\,|x|^{A-\sigma}dx
\end{equation}
for all $0\le\varphi\in H^1_c(G)\cap L^\infty_c(G)$.
The notions of a nonnegative sub-solution and solution are defined similarly
by replacing "$\ge$" with "$\le$" and "$=$" respectively.
By the weak Harnack inequality for super-solutions
any nontrivial nonnegative super-solution $u$ to \eqref{e:MAIN} in $G$ is strictly positive in $G$,
in the sense that $u^{-1}\in L^\infty_{loc}(G)$.
\smallskip

%Equation \eqref{e:MAIN} is equivalent to a semilinear equation
%\begin{equation}\label{e:DRIFT}
%-\Delta u-\frac{Ax}{|x|^2}\cdot\nabla u-\frac{B}{|x|^2}u=\frac{c}{|x|^\sigma}u^p
%\quad\mbox{in }\:\C^\rho_\Omega.
%\end{equation}
%A {\em (super)\,solution} to \eqref{e:DRIFT} is an $u\in H^1_{loc}(G)$ such that
%\begin{equation}\label{e:DefSol-DRIFT}
%\int_G\nabla u\cdot\nabla\varphi\,dx-\int_G\nabla
%u\cdot\frac{Ax}{|x|^2}\varphi\,dx- \int_G\frac{B}{|x|^2}u\varphi\,dx\,
%(\geq)=\,\int_G\frac{c}{|x|^\sigma}u^p\varphi\,dx,
%\end{equation}
%for all $0\le\varphi\in H^1_c(\C_\Omega^\rho)\cap L^\infty_c(\C_\Omega^\rho)$.

Equation \eqref{e:MAIN} comprising in
particular the known in astrophysics Lane-Emdem equation, is a
prototype  model for general semilinear equations. The qualitative
theory of equations of type \eqref{e:MAIN} has been extensively
studied because of their rich mathematical structure  and various
applications for the whole range of the parameter $p\in\R$, e.g.
in combustion theory ($p>1$) \cite{Galaktionov}, population dynamics ($0<p<1$)
\cite{Murray}, pseudoplastic fluids ($p<0$) \cite{Hernandez,Lazer}.
It has been known at least since earlier works by Serrin (cf. the
references in \cite{Serrin-Zou}) and celebrated paper by Gidas and
Spruck \cite{Gidas-Spruck} that equations of type \eqref{e:MAIN}
on unbounded domains admit positive (super)\,solutions only for specific values
of $(p,\sigma)\in\R^2$.
%The equation of type \eqref{e:MAIN} have been studied so far  mainly for super-linear case $p>1$.
%(see, e.g. \cite{Bandle-Levine,Brezis-Cabre,Brezis-Kamin,Dupaigne,Gidas-Spruck,Levine,Smets,Pohozaev-Tesei}
%and references therein).
For instance, it is well known by now that the equation
\begin{equation}\label{e:KLMS}
-\Delta u=u^p
\end{equation}
in exterior of a ball in $\R^N$ ($N\geq 3$) has no positive
super-solutions if $p\leq\frac{N}{N-2}$. The {\it critical exponent} $p^*=\frac{N}{N-2}$
is sharp in the sense that it separates the zones of the existence
and nonexistence, i.e. for $p>p^*$ \eqref{e:MAIN} has positive solutions outside a ball.
This result has been extended in several directions
(see, e.g. \cite{Bandle-Levine,BCN,B-VP,Birindelli,Brezis-Cabre,
Brezis-Kamin,Dupaigne,KLS,KLSU,Levine,Serrin-Zou,Smets,Pohozaev-Tesei,Zhang-1,Zhang-2},
references therein, and the list is by no means complete).
In particular, in \cite{KLS} it was
shown that the critical exponent $p^\ast=\frac{N}{N-2}$ is stable
with respect to the change of the Laplacian by a second--order
uniformly elliptic divergence type operator with measurable
coefficients $-\sum\partial_i(a_{ij}\partial_j)$, perturbed by a potential,
for a sufficiently wide class of potentials.
For instance, for $\epsilon>0$ the equation
\begin{equation}\label{e:EPS}
-\Delta u-\frac{B}{|x|^{2+\epsilon}}u=u^p
\end{equation}
in the exterior of a ball in $\R^N$ ($N\geq 3$)
has the same critical exponent as \eqref{e:KLMS} \cite[Theorem 1.2]{KLS}.
On the other hand it is easy to see that if $\epsilon<0$ and $B>0$
then \eqref{e:EPS} has no positive super-solutions for any
$p\in\R$, while if $\epsilon<0$ and $B<0$ then \eqref{e:EPS} admits positive solutions
for all $p\in\R$ ($p\neq 1$).
In the borderline case $\epsilon=0$ the critical exponent $p^\ast$ becomes explicitly dependent
on the parameter $B$. This phenomenon and its relation with Hardy type inequalities has
been recently observed on a ball and/or exterior domains in \cite{Brezis-Cabre,Dupaigne,Terracini}
in the case $p>1$.

The equation with first order term
\begin{equation}\label{e:drift}
-\Delta{u}-\frac{Ax}{|x|^{2+\epsilon}}\nabla{u}=u^p
\end{equation}
in the exterior of a ball in $\R^N$ ($N\geq 3$) represents another type of behavior.
If $\epsilon>0$ then \eqref{e:drift} has the same critical exponent $p^\ast=\frac{N}{N-2}$
as \eqref{e:KLMS}, and $p^\ast$ is stable with respect to the change of the Laplacian
by a second--order uniformly elliptic divergence type operator \cite[Theorem 1.8]{KLSU}.
On the other hand it is easy to see that if $\epsilon<0$ and $A>0$ then
\eqref{e:drift} has no positive super-solutions if and only if $p\le 1$,
while if $\epsilon<0$ and $A<0$ then \eqref{e:drift} has no positive super-solutions
if and only if $p\ge 1$.
In the borderline case $\epsilon=0$ the critical exponent $p^\ast$ explicitly depends
on the parameter $A$ (see \cite{Pohozaev-Tesei,Smets} for the case $p>1$).

When considered on cone--like domains, the nonexistence zone
depends in addition on the principal Dirichlet eigenvalue of the cross-section of the cone.
In the super-linear case $p>1$ the equation
\begin{equation}\label{e:cone}
-\Delta{u}=u^p\quad \text{in }\:\C_\Omega^1
\end{equation}
has been considered in \cite{Bandle,Bandle-Levine,BCN}
(see also \cite{Birindelli} for systems and \cite{KLM} for uniformly elliptic equations with
measurable coefficients). A new nonexistence phenomenon for the
sublinear case $p<1$ has been recently revealed in \cite{KLMS}.
Particularly, it was discovered that equation \eqref{e:cone}
in a proper cone-like domain has two critical exponents,
the second one appearing in the sublinear case, so that \eqref{e:cone}
has no positive super-solutions if and only if $p_\ast\le p\le p^\ast$,
where $p_\ast<1$ and $p^\ast>1$.
In \cite{KLM} for $p>1$ it was shown that if the Laplacian is replaced
by a second--order uniformly elliptic divergence type operator $-\sum\partial_i(a_{ij}\partial_j)$
then the value of the critical exponents on the cone depends on the coefficients of the matrix
$(a_{ij}(x))$ as well as on the geometry of the cross--section.

In the present paper we study equation \eqref{e:MAIN} on cone-like domains for the full range of
the parameters $p,\sigma,A,B$.
Note that \eqref{e:MAIN} can be rewritten in the form
$$-\Delta{u}-\frac{Ax}{|x|^2}\nabla{u}-\frac{B}{|x|^2}{u}=
\frac{C}{|x|^\sigma}{u}^p\quad \text{in }\quad{\C}_\Omega^{\rho},$$
so it represents the borderline case both with respect to the zero order
and the first order perturbations in the linear part.
As we will see below, due to the presence of the weighted function
and lower order terms equation \eqref{e:MAIN} exhibits all the cases of qualitative
behavior described above for the Laplacian. Our approach
to the problem in this paper is the development of the method introduced
in \cite{KLS} (see also \cite{KLS,KLSU,KLM,KLMS}) and is different from the techniques used in
\cite{Bandle,BCN,B-VP,Brezis-Cabre,Dupaigne,Pohozaev-Tesei,Smets}.
It is based on the explicit construction of appropriate barriers and involves
the analysis of asymptotic behavior of super-harmonic functions associated to the
Laplace operator with critical potentials, Phragmen--Lindel\"of type comparison arguments
and an improved version of Hardy's inequality in cone--like domains.
The advantages of our approach are its transparency and
flexibility. Particularly we prove the nonexistence results for
the most general definition of weak solutions and avoid any
assumptions on the smoothness of the boundary of the cone.

Below we denote $C_H:=\frac{(2-N-A)^2}{4}$, while
$\lambda_1=\lambda_1(\Omega)\ge 0$ denotes the principal Dirichlet
eigenvalue of the Laplace--Beltrami operator $-\Delta_\omega$ on $\Omega$.
First, we formulate the result in the special linear case.
\begin{theorem}\label{t:MAIN-0}
Let $(p,\sigma)=(1,2)$.
Then equation \eqref{e:MAIN} has no positive super-solutions if and only if $B+C>C_H+\lambda_1$.
\end{theorem}
If $B\le C_H+\lambda_1$ then the quadratic equation
\begin{equation}\label{e:ROOTS2}
\gamma(\gamma+N-2+A)=\lambda_1-B
\end{equation}
has real roots, denoted by $\gamma^-\le\gamma^+$.
Note that if $B=C_H+\lambda_1$, then $\gamma^\pm=(2-N-A)/2$.
%By a standard argument (cf. Corollary \ref{APA})
%inequality \eqref{e:HARDY} implies that if \eqref{e:ROOTS2} has no real roots then
%\eqref{e:SING} has no positive super-solutions for any $(p,\sigma)\in\R^2$.
%When $(p,\sigma)=(1,2)$ then equation \eqref{e:SING} has no positive super-solutions
%if and only if $c\geq \lambda_1+c_H-B$, where $c>0$ is a constant in \eqref{e:MAIN}.
%We confine ourself to the case when $\gamma^-,\gamma^+\in\R$ and $(p,\sigma)\in\R^2_\ast$,
%where $\R^2_\ast:=\R^2\setminus\{1,2\}$.

For $B\le C_H+\lambda_1$ we introduce the critical line
$\Lambda_\ast(p,A,B,\Omega)$ on the $(p,\sigma)$--plane
$$\Lambda_\ast(p,A,B,\Omega):=\min\{\gamma^-(p-1)+2,\;\gamma^+(p-1)+2\}\qquad (p\in\R),$$
and the nonexistence set
$$\mathcal{N}=\{(p,\sigma)\in\R^2\setminus\{1,2\}:\textit{equation \eqref{e:MAIN} has no positive super-solutions}\}.$$
The main result of the paper reads as follows.

\begin{theorem}\label{t:MAIN}
The following assertions are valid.
\begin{enumerate}
\item[{$(i)$}]
Let $B<C_H+\lambda_1$. Then $\mathcal{N}=\{\sigma\leq \Lambda_\ast(p)\}$.
\item[{$(ii)$}]
Let $B=C_H+\lambda_1$. Then
$$\{\sigma<\Lambda_\ast(p)\}\cup\{\sigma=\Lambda_\ast(p),\,p\geq -1\}\subseteq\mathcal{N}\subseteq\{\sigma\le\Lambda_\ast(p)\}.$$
If $\,\Omega=S^{N-1}$ then
$\mathcal{N}=\{\sigma<\Lambda_\ast(p)\}\cup\{\sigma=\Lambda_\ast(p),\,p\geq
-1\}$.
\end{enumerate}
\end{theorem}

%%%%%%%%%%%%%%%%% Begin of eepic-graphic %%%%%%%%%%%%%%%%%%%%%

\setlength{\unitlength}{0.0004in}
\begingroup\makeatletter\ifx\SetFigFont\undefined%
\gdef\SetFigFont#1#2#3#4#5{%
  \reset@font\fontsize{#1}{#2pt}%
  \fontfamily{#3}\fontseries{#4}\fontshape{#5}%
  \selectfont}%
\fi\endgroup%
\renewcommand{\dashlinestretch}{30}

%%%%%%%%%%%%%%%%%%%%%%%%%%%%%%%%%%%%%%%%%%%%%%%%%%%%

\begin{figure}[p]
\begin{center}

\vspace{1cm}

%
%%%%% Row 1
%
%

\begin{picture}(4815,4800)(0,-300) %% $\gamma_-<0$, $\gamma_+\ge 0$
\put(1000,-300){\makebox(0,0)[lb]{\smash{{{\SetFigFont{9}{12.0}{\rmdefault}{\mddefault}{\updefault}$(a)\;$: $\;\gamma_-<0$, $\gamma_+\ge 0$}}}}}

\texture{88555555 55000000 555555 55000000 555555 55000000 555555 55000000
           555555 55000000 555555 55000000 555555 55000000 555555 55000000
           555555 55000000 555555 55000000 555555 55000000 555555 55000000
           555555 55000000 555555 55000000 555555 55000000 555555 55000000 }
\thinlines\shade\path(3012,3012)(4512,12)(12,12)(12,612)(3012,3012)

\path(12,1812)(4812,1812) % OP-axis
\blacken\path(4692.000,1782.000)(4812.000,1812.000)(4692.000,1842.000)(4692.000,1782.000) % arrow
\path(2412,12)(2412,4762) % OS-axis
\blacken\path(2442.000,4642.000)(2412.000,4762.000)(2382.000,4642.000)(2442.000,4642.000) % arrow

\dottedline{120}(12,3012)(4812,3012)
\dottedline{120}(3012,4762)(3012,12)

\Thicklines
\drawline(12,612)(3012,3012)(4512,12) % bold line nonexistence area border
\put(3012,3012){\whiten\ellipse{100}{100}}
\put(1510,1812){\ellipse{50}{50}}
\put(3620,1812){\ellipse{50}{50}}
\put(2412,2540){\ellipse{50}{50}}

\put(4690,1630){\makebox(0,0)[lb]{\smash{{{\SetFigFont{7}{12.0}{\rmdefault}{\mddefault}{\updefault}$p$}}}}}
\put(2180,4650){\makebox(0,0)[lb]{\smash{{{\SetFigFont{7}{12.0}{\rmdefault}{\mddefault}{\updefault}$\sigma$}}}}}
\put(2200,3070){\makebox(0,0)[lb]{\smash{{{\SetFigFont{7}{12.0}{\rmdefault}{\mddefault}{\updefault}$2$}}}}}
\put(2840,1590){\makebox(0,0)[lb]{\smash{{{\SetFigFont{7}{12.0}{\rmdefault}{\mddefault}{\updefault}$1$}}}}}
\put(3580,2030){\makebox(0,0)[lb]{\smash{{{\SetFigFont{7}{12.0}{\rmdefault}{\mddefault}{\updefault}$1\!-\!\frac{2}{\gamma_-}$}}}}}
\put(770,2030){\makebox(0,0)[lb]{\smash{{{\SetFigFont{7}{12.0}{\rmdefault}{\mddefault}{\updefault}$1\!-\!\frac{2}{\gamma_+}$}}}}}
\put(1650,2600){\makebox(0,0)[lb]{\smash{{{\SetFigFont{7}{12.0}{\rmdefault}{\mddefault}{\updefault}$2\!-\!\gamma_+$}}}}}
\put(300,300){\makebox(0,0)[lb]{\smash{{{\SetFigFont{9}{12.0}{\rmdefault}{\mddefault}{\itdefault}$\mathcal{N}$}}}}}
\end{picture}
\hspace{3cm}
\begin{picture}(4815,4800)(0,-300) %%% $\gamma_-=\gamma_+=0$
\put(1000,-300){\makebox(0,0)[lb]{\smash{{{\SetFigFont{9}{12.0}{\rmdefault}{\mddefault}{\updefault}$(b)\;$: $\;\gamma_-=\gamma_+=0$}}}}}

\texture{88555555 55000000 555555 55000000 555555 55000000 555555 55000000
           555555 55000000 555555 55000000 555555 55000000 555555 55000000
           555555 55000000 555555 55000000 555555 55000000 555555 55000000
           555555 55000000 555555 55000000 555555 55000000 555555 55000000 }
\thinlines\shade\path(3012,3012)(4812,3012)(4812,12)(12,12)(12,3012)(3012,3012)

\path(12,1812)(4812,1812) % OP-axis
\blacken\path(4692.000,1782.000)(4812.000,1812.000)(4692.000,1842.000)(4692.000,1782.000) % arrow
\path(2412,12)(2412,4762) % OS-axis
\blacken\path(2442.000,4642.000)(2412.000,4762.000)(2382.000,4642.000)(2442.000,4642.000) % arrow

\dottedline{120}(3012,4762)(3012,12)
\dottedline{60}(1850,3012)(1850,1812)

\Thicklines
\drawline(1850,3012)(4812,3012) % bold line nonexistence area border
\Thicklines
\dashline[50]{200}(12,3012)(1600,3012) % dash line nonexistence area border
%%\dottedline{180}(12,3012)(1750,3012) % dotted line nonexistence area border
\put(1850,3012){\ellipse{60}{60}}

\put(4600,1600){\makebox(0,0)[lb]{\smash{{{\SetFigFont{7}{12.0}{\rmdefault}{\mddefault}{\updefault}$p$}}}}}
\put(2180,4650){\makebox(0,0)[lb]{\smash{{{\SetFigFont{7}{12.0}{\rmdefault}{\mddefault}{\updefault}$\sigma$}}}}}
\put(2200,3070){\makebox(0,0)[lb]{\smash{{{\SetFigFont{7}{12.0}{\rmdefault}{\mddefault}{\updefault}$2$}}}}}
\put(2840,1590){\makebox(0,0)[lb]{\smash{{{\SetFigFont{7}{12.0}{\rmdefault}{\mddefault}{\updefault}$1$}}}}}
\put(1650,1590){\makebox(0,0)[lb]{\smash{{{\SetFigFont{7}{12.0}{\rmdefault}{\mddefault}{\updefault}$-\!1$}}}}}
\put(300,300){\makebox(0,0)[lb]{\smash{{{\SetFigFont{9}{12.0}{\rmdefault}{\mddefault}{\itdefault}$\mathcal{N}$}}}}}
\end{picture}

\vspace{2cm}
%
%
%%%%% Row 2
%
%
\begin{picture}(4815,4800)(0,-300) %% $\gamma_-,\gamma_+<0$
\put(1200,-300){\makebox(0,0)[lb]{\smash{{{\SetFigFont{9}{12.0}{\rmdefault}{\mddefault}{\updefault}$(e):\;$ $\;\gamma_-,\gamma_+<0$}}}}}

\texture{88555555 55000000 555555 55000000 555555 55000000 555555 55000000
           555555 55000000 555555 55000000 555555 55000000 555555 55000000
           555555 55000000 555555 55000000 555555 55000000 555555 55000000
           555555 55000000 555555 55000000 555555 55000000 555555 55000000 }
\thinlines\shade\path(3012,3012)(4512,12)(12,12)(12,4762)(1012,4762)(3012,3012)

\path(12,1812)(4812,1812) % OP-axis
\blacken\path(4692.000,1782.000)(4812.000,1812.000)(4692.000,1842.000)(4692.000,1782.000) % arrow
\path(2412,12)(2412,4762) % OS-axis
\blacken\path(2442.000,4642.000)(2412.000,4762.000)(2382.000,4642.000)(2442.000,4642.000) % arrow

\dottedline{120}(12,3012)(4812,3012)
\dottedline{120}(3012,4762)(3012,12)

\Thicklines
\drawline(1012,4762)(3012,3012)(4512,12)
\put(3012,3012){\whiten\ellipse{100}{100}}
\put(3620,1812){\ellipse{50}{50}}
\put(2412,3532){\ellipse{50}{50}}

\put(4690,1630){\makebox(0,0)[lb]{\smash{{{\SetFigFont{7}{12.0}{\rmdefault}{\mddefault}{\updefault}$p$}}}}}
\put(2180,4650){\makebox(0,0)[lb]{\smash{{{\SetFigFont{7}{12.0}{\rmdefault}{\mddefault}{\updefault}$\sigma$}}}}}
\put(2200,3070){\makebox(0,0)[lb]{\smash{{{\SetFigFont{7}{12.0}{\rmdefault}{\mddefault}{\updefault}$2$}}}}}
\put(2840,1590){\makebox(0,0)[lb]{\smash{{{\SetFigFont{7}{12.0}{\rmdefault}{\mddefault}{\updefault}$1$}}}}}
\put(3580,2030){\makebox(0,0)[lb]{\smash{{{\SetFigFont{7}{12.0}{\rmdefault}{\mddefault}{\updefault}$1\!-\!\frac{2}{\gamma_-}$}}}}}
\put(1630,3480){\makebox(0,0)[lb]{\smash{{{\SetFigFont{7}{12.0}{\rmdefault}{\mddefault}{\updefault}$2\!-\!\gamma_+$}}}}}
\put(300,300){\makebox(0,0)[lb]{\smash{{{\SetFigFont{9}{12.0}{\rmdefault}{\mddefault}{\itdefault}$\mathcal{N}$}}}}}
\end{picture}
\hspace{3cm}
\begin{picture}(4815,4800)(0,-300) %% $\gamma_-=\gamma_+<0$
\put(1000,-300){\makebox(0,0)[lb]{\smash{{{\SetFigFont{9}{12.0}{\rmdefault}{\mddefault}{\updefault}$(f):\;$ $\;\gamma_-=\gamma_+<0$}}}}}
\texture{88555555 55000000 555555 55000000 555555 55000000 555555 55000000
           555555 55000000 555555 55000000 555555 55000000 555555 55000000
           555555 55000000 555555 55000000 555555 55000000 555555 55000000
           555555 55000000 555555 55000000 555555 55000000 555555 55000000 }
\thinlines\shade\path(12,12)(12,4762)(1012,4762)(4812,1312)(4812,12)(12,12)

\path(12,1812)(4812,1812) % OP-axis
\blacken\path(4692.000,1782.000)(4812.000,1812.000)(4692.000,1842.000)(4692.000,1782.000) % arrow
\path(2412,12)(2412,4762) % OS-axis
\blacken\path(2442.000,4642.000)(2412.000,4762.000)(2382.000,4642.000)(2442.000,4642.000) % arrow

\dottedline{120}(12,3012)(4812,3012)
\dottedline{120}(3012,4762)(3012,12)
\dottedline{60}(1812,4042)(1812,1812)
\dottedline{60}(1812,4042)(2412,4042)

\Thicklines
\drawline(1812,4042)(4812,1312) % bold line nonexistence area border
\Thicklines
\dashline[50]{200}(1012,4762)(1650,4180) % dashline nonexistence area border
%%\dottedline{180}(1012,4762)(1650,4180)
\put(1812,4042){\ellipse{60}{60}}

\put(4690,1630){\makebox(0,0)[lb]{\smash{{{\SetFigFont{7}{12.0}{\rmdefault}{\mddefault}{\updefault}$p$}}}}}
\put(2180,4650){\makebox(0,0)[lb]{\smash{{{\SetFigFont{7}{12.0}{\rmdefault}{\mddefault}{\updefault}$\sigma$}}}}}
\put(2200,3070){\makebox(0,0)[lb]{\smash{{{\SetFigFont{7}{12.0}{\rmdefault}{\mddefault}{\updefault}$2$}}}}}
\put(2440,4000){\makebox(0,0)[lb]{\smash{{{\SetFigFont{6}{12.0}{\rmdefault}{\mddefault}{\updefault}$A\!\!+\!\!N$}}}}}
\put(2840,1590){\makebox(0,0)[lb]{\smash{{{\SetFigFont{7}{12.0}{\rmdefault}{\mddefault}{\updefault}$1$}}}}}
\put(1650,1590){\makebox(0,0)[lb]{\smash{{{\SetFigFont{7}{12.0}{\rmdefault}{\mddefault}{\updefault}$-\!1$}}}}}
\put(300,300){\makebox(0,0)[lb]{\smash{{{\SetFigFont{9}{12.0}{\rmdefault}{\mddefault}{\itdefault}$\mathcal{N}$}}}}}
\end{picture}

\vspace{2cm}
%
%
%%%%% Row 3
%
%
\begin{picture}(4815,4800)(0,-300) %% $\gamma_-,\gamma_+>0$
\put(1000,-300){\makebox(0,0)[lb]{\smash{{{\SetFigFont{9}{12.0}{\rmdefault}{\mddefault}{\updefault}$(c):\;$ $\;\gamma_-\ge 0$, $\gamma_+>0$}}}}}

\texture{88555555 55000000 555555 55000000 555555 55000000 555555 55000000
           555555 55000000 555555 55000000 555555 55000000 555555 55000000
           555555 55000000 555555 55000000 555555 55000000 555555 55000000
           555555 55000000 555555 55000000 555555 55000000 555555 55000000 }
\thinlines\shade\path(3012,3012)(4812,3412)(4812,12)(12,12)(12,612)(3012,3012)

\path(12,1812)(4812,1812) % OP-axis
\blacken\path(4692.000,1782.000)(4812.000,1812.000)(4692.000,1842.000)(4692.000,1782.000) % arrow
\path(2412,12)(2412,4762) % OS-axis
\blacken\path(2442.000,4642.000)(2412.000,4762.000)(2382.000,4642.000)(2442.000,4642.000) % arrow

\dottedline{120}(12,3012)(4650,3012)
\dottedline{120}(3012,4762)(3012,12)

\Thicklines
\path(12,612)(3012,3012)(4812,3412)
\put(3012,3012){\whiten\ellipse{100}{100}}
\put(1510,1812){\ellipse{50}{50}}
\put(2412,2540){\ellipse{50}{50}}

\put(4600,1600){\makebox(0,0)[lb]{\smash{{{\SetFigFont{7}{12.0}{\rmdefault}{\mddefault}{\updefault}$p$}}}}}
\put(2180,4650){\makebox(0,0)[lb]{\smash{{{\SetFigFont{7}{12.0}{\rmdefault}{\mddefault}{\updefault}$\sigma$}}}}}
\put(2200,3070){\makebox(0,0)[lb]{\smash{{{\SetFigFont{7}{12.0}{\rmdefault}{\mddefault}{\updefault}$2$}}}}}
\put(2840,1590){\makebox(0,0)[lb]{\smash{{{\SetFigFont{7}{12.0}{\rmdefault}{\mddefault}{\updefault}$1$}}}}}
\put(770,2030){\makebox(0,0)[lb]{\smash{{{\SetFigFont{7}{12.0}{\rmdefault}{\mddefault}{\updefault}$1\!-\!\frac{2}{\gamma_+}$}}}}}
\put(1650,2600){\makebox(0,0)[lb]{\smash{{{\SetFigFont{7}{12.0}{\rmdefault}{\mddefault}{\updefault}$2\!-\!\gamma_+$}}}}}
\put(300,300){\makebox(0,0)[lb]{\smash{{{\SetFigFont{9}{12.0}{\rmdefault}{\mddefault}{\itdefault}$\mathcal{N}$}}}}}
\end{picture}
\hspace{3cm}
\begin{picture}(4815,4800)(0,-300) %% $\gamma_-=\gamma_+>0$
\put(1000,-300){\makebox(0,0)[lb]{\smash{{{\SetFigFont{9}{12.0}{\rmdefault}{\mddefault}{\updefault}$(d):\;$ $\;\gamma_-=\gamma_+>0$}}}}}

\texture{88555555 55000000 555555 55000000 555555 55000000 555555 55000000
           555555 55000000 555555 55000000 555555 55000000 555555 55000000
           555555 55000000 555555 55000000 555555 55000000 555555 55000000
           555555 55000000 555555 55000000 555555 55000000 555555 55000000 }
\thinlines\shade\path(3012,3012)(4812,4046)(4812,12)(12,12)(12,1212)(3012,3012)

\path(12,1812)(4812,1812) % OP-axis
\blacken\path(4692.000,1782.000)(4812.000,1812.000)(4692.000,1842.000)(4692.000,1782.000) % arrow
\path(2412,12)(2412,4762) % OS-axis
\blacken\path(2442.000,4642.000)(2412.000,4762.000)(2382.000,4642.000)(2442.000,4642.000) % arrow

\dottedline{120}(12,3012)(4812,3012)
\dottedline{120}(3012,4762)(3012,12)
\dottedline{60}(1850,2310)(1850,1812)
\dottedline{60}(1850,2310)(2412,2310)

\Thicklines
\drawline(1850,2315)(4812,4046) % bold line nonexistence area border
\Thicklines
\dashline[50]{200}(12,1212)(1650,2190) % dashline nonexistence area border
%%\dottedline{180}(12,1212)(1650,2190)
\put(1850,2310){\ellipse{60}{60}}

\put(4600,1630){\makebox(0,0)[lb]{\smash{{{\SetFigFont{7}{12.0}{\rmdefault}{\mddefault}{\updefault}$p$}}}}}
\put(2180,4650){\makebox(0,0)[lb]{\smash{{{\SetFigFont{7}{12.0}{\rmdefault}{\mddefault}{\updefault}$\sigma$}}}}}
\put(2200,3070){\makebox(0,0)[lb]{\smash{{{\SetFigFont{7}{12.0}{\rmdefault}{\mddefault}{\updefault}$2$}}}}}
\put(2440,2270){\makebox(0,0)[lb]{\smash{{{\SetFigFont{6}{12.0}{\rmdefault}{\mddefault}{\updefault}$A\!\!+\!\!N$}}}}}
\put(2840,1590){\makebox(0,0)[lb]{\smash{{{\SetFigFont{7}{12.0}{\rmdefault}{\mddefault}{\updefault}$1$}}}}}
\put(1650,1590){\makebox(0,0)[lb]{\smash{{{\SetFigFont{7}{12.0}{\rmdefault}{\mddefault}{\updefault}$-\!1$}}}}}
\put(300,300){\makebox(0,0)[lb]{\smash{{{\SetFigFont{9}{12.0}{\rmdefault}{\mddefault}{\itdefault}$\mathcal{N}$}}}}}

\end{picture}
\vspace{1cm}

\caption{The nonexistence set ${\mathcal N}$ of equation
\eqref{e:MAIN} for typical values of $\gamma^-$ and
$\gamma^+$.}\label{fig}
\end{center}
\end{figure}

%%%%%%%%%%%%%%%%% End of graphic %%%%%%%%%%%%%%%%%%%%%

\begin{remark}
$(i)$ Observe that the nonexistence set ${\mathcal N}$ does not depend on the
value of the parameter $C>0$ in \eqref{e:MAIN}. In view of the
scaling invariance of \eqref{e:MAIN}
the set ${\mathcal N}$ also does not depend on the value of $\rho>0$.\\
$(ii)$ Using sub and super-solutions techniques one can show that
if \eqref{e:MAIN} has a positive super-solution in
$\C_\Omega^\rho$ then it has a positive solution in
$\C_\Omega^\rho$ (cf. \cite[Proposition 1.1]{KLM}). Thus for any
$(p,\sigma)\in\R^2\setminus{\mathcal N}$
equation \eqref{e:MAIN} admits positive solutions.\\
$(iii)$ In the case of proper domains $\Omega\Subset S^{N-1}$, the
existence (or nonexistence) of positive super-solutions to
\eqref{e:main'} with $p<-1$ and $s=\alpha_\ast(p-1)+2$ becomes a
more involved issue that remains open at the moment.
%The analysis of the decay rate of super-solutions near the lateral boundary of the cone
%should be invoked (see discussion at the end of Section~6).
We will return to this problem elsewhere.
\end{remark}

\begin{remark}
Figure \ref{fig} shows the qualitative pictures of the set
${\mathcal N}$ for typical values of $\gamma^-$, $\gamma^+$. The
case $(a)$ is typical for $A,B=0$. The case $(b)$ occurs, e.g.,
when $A,B=0$ and $N=2$.
The cases $(c)$ and $(d)$ appear, in particular, when $A=0$ ($B<C_H$ and $B=C_H$ respectively).
The cases $(e)$ and $(f)$ are never realized by \eqref{e:MAIN} with $A=0$.
Assume, for instance, that $B=0$, $\lambda_1=0$ and $\sigma=0$.
Then \eqref{e:MAIN} admits at most one critical exponent $p^\ast$ which,
depending whether $N+A>2$ or $N+A<2$, appears in the superlinear case ($p>1$) or sublinear case ($p<1$).
In the former case there are no positive super-solutions if and only if $p\le p^*$,
whereas in the latter if and only if $p^*\le p$.
Thus $N+A$ plays a role of the "effective dimension".
Similar behavior is exhibited by second-order elliptic nondivergent type equations
with measurable coefficients
$-\sum a_{ij}\partial^2_{ij} u=u^p$
in the exterior of a ball in $\R^N$, which were recently studied in \cite{KLS1}.
The value of the critical exponent for such equations depends on the behavior
of the matrix $(a_{ij}(x))$ at infinity, though not directly but
via an "effective dimension" which is determined by the asymptotic of $(a_{ij}(x))$.
%as $$\Psi=\frac{Tr(a)}{\frac{(ax,x)}{|x|^2}}.$$
\end{remark}

Applying the Kelvin transformation $y=y(x)=\frac{x}{|x|^2}$
one sees that if $u$ is a positive (super)\,solution to equation (\ref{e:MAIN})
then $\check{u}(y)=|y|^{2-N}u(x(y))$ is a positive (super)\,solution to the equation
\begin{equation}\label{e:KELVIN}
-\nabla\cdot(|x|^A\nabla
\check{u})-B|x|^{A-2}\check{u}=C|x|^{A-s}\check{u}^p\quad \text{in
}\quad\check{\C}_\Omega^{1},
\end{equation}
where $s=(N+2)-p(N-2)-\sigma$ and
$\check{\C}_\Omega^{1}:=\{(r,\omega)\in\R^N:\,\omega\in\Omega,\,0<r<1\}$
is the interior cone--like domain.
For equation \eqref{e:KELVIN} we define the critical line
$$\Lambda^*(p,A,B,\Omega):=\max\{\gamma^-(p-1)+2,\;\gamma^+(p-1)+2\}\qquad (p\in\R),$$
and the set
$\check{{\mathcal N}}=\{(p,s)\in\R^2\setminus\{1,2\}:\textit{\eqref{e:KELVIN} has no positive super-solutions}\}$.
The following theorem
extends the results in \cite{Terracini,Brezis-Cabre,Dupaigne} ($A=0$)
and \cite{Smets,Pohozaev-Tesei} ($B=0$), obtained on the punctured ball
in the super-linear case $p>1$. It is derived from Theorem \ref{t:MAIN} via the Kelvin transformation.

\begin{theorem}\label{t:KELVIN}
The following assertions are valid.
\begin{enumerate}
\item[{$(i)$}]
Let $B<C_H+\lambda_1$. Then $\mathcal{N}=\{s\geq
\Lambda^\ast(p)\}$.
\item[{$(ii)$}]
Let $B=C_H+\lambda_1$. Then
$$\{s>\Lambda^\ast(p)\}\cup\{s=\Lambda^\ast(p),\,p\geq -1\}\subseteq\mathcal{N}\subseteq\{s\ge\Lambda^\ast(p)\}.$$
If $\,\Omega=S^{N-1}$ then
$\mathcal{N}=\{s>\Lambda^\ast(p)\}\cup\{s=\Lambda^\ast(p),\,p\geq
-1\}$
\end{enumerate}
\end{theorem}

%\begin{remark}
%Theorem \ref{t:KELVIN} extends the results in \cite{Terracini,Brezis-Cabre,Dupaigne} ($A=0$)
%and \cite{Smets,Pohozaev-Tesei} ($B=0$), obtained on the punctured ball
%in the super-linear case $p>1$.
%\end{remark}

Theorem \ref{t:MAIN} is proved in the paper after a reduction of \eqref{e:MAIN}
to the uniformly elliptic case $A=0$.
This reduction is described in Section 2 below.
The rest of the paper is organized  as follows.
In Section~3 we prove a version of the improved Hardy inequality in cone--like domains.
In Section~4 we study asymptotical behavior of super-solutions to certain linear equations.
The proof of the main result is contained in Section~5 (super-linear case $p\ge 1$)
and Section~6 (sub-linear case $p<1$).
Finally, Appendix includes auxiliary results on the relation between
the existence of positive solutions to linear equations and
positivity properties of the corresponding quadratic forms.

\section{Equivalent statement of the problem}

The next lemma shows that a simple transformation allows one to reduce equation
(\ref{e:MAIN}) to the uniformly elliptic case $A=0$.
%Direct computation verifies the proof of the following lemma.
\begin{lemma} The function $u$ is a (super)\,solutions to equation \eqref{e:MAIN}
if and only if $w(x)=|x|^{\frac{A}{2}}u(x)$ is a (super)\,solution
to the equation
\begin{equation}\label{e:main'}
-\Delta w-\frac{\mu}{|x|^2}w=\frac{C}{|x|^s}w^p\quad\text{in
}\:\C_\Omega^\rho,
\end{equation}
where $\mu=B-\frac{A}{2}(\frac{A}{2}+N-2)$ and
$s=\sigma+\frac{A}{2}(p-1)$.
%The (super)\,solutions to equation \eqref{e:MAIN} are in
%one--to--one correspondence to (super)\, solutions to the equation
%\begin{equation}\label{e:main'}
%-\Delta w-\frac{\mu}{|x|^2}w=\frac{c}{|x|^s}w^p\quad\text{in
%}\:\C_\Omega^\rho,
%\end{equation}
%where $\mu=B-\frac{A}{2}(\frac{A}{2}+N-2)$ and
%$s=\sigma+\frac{A}{2}(p-1)$. The correspondence is given by
%$$u(x)=|x|^{-\frac{A}{2}}w(x).$$
\end{lemma}

\begin{proof}
The direct computation.
\end{proof}

The existence of positive solutions to \eqref{e:main'} is intimately related
to an associated Hardy type inequality for exterior cone--like domains,
which has the form
\begin{equation}\label{e:HARDY}
\int_{\C_\Omega^\rho}|\nabla u|^2dx\geq
(C_H+\lambda_1)\int_{\C_\Omega^\rho}\frac{u^2}{|x|^{2}}dx
%+c_\ast\int_{\C_\Omega^\rho}\frac{u^2}{\log^2|x|}|x|^{A-2}dx
,\qquad\forall u\in C^\infty_c(\C_\Omega^\rho),
\end{equation}
where $C_H:=\frac{(N-2)^2}{4}$ and the constant $C_H+\lambda_1$ is
sharp. We prove a refined version of \eqref{e:HARDY} in Section~\ref{s:Hardy}.
By virtue of Lemma \ref{l:A-ground} in Appendix, inequality \eqref{e:HARDY} implies
that equation \eqref{e:main'} has positive super-solutions if and only if
$\mu\le C_H+\lambda_1$, see Remark \ref{r:Nonexist} below.
Theorem \ref{t:MAIN-0} is an immediate consequence of this result.

If $\mu\le C_H+\lambda_1$ then the quadratic equation
\begin{equation}\label{e:roots'}
\alpha(\alpha+N-2)=\lambda_1-\mu.
\end{equation}
has real roots, denoted by $\alpha^-\le\alpha^+$. If
$\mu=C_H+\lambda_1$ then $\alpha^\pm=\frac{2-N}{2}$.
In this case we write $\alpha_\ast:=\frac{2-N}{2}$ for convenience. As before, we
introduce the critical line
$$\Lambda=\Lambda(p,\mu,\Omega):=\min\{\alpha^-(p-1)+2,\;\alpha^+(p-1)+2\}\qquad (p\in\R),$$
and the nonexistence set
$$\mathcal{N}=\{(p,s)\in\R^2\setminus\{1,2\}:\textit{equation \eqref{e:main'} has no positive super-solutions}\}.$$
Theorem \ref{t:MAIN} is a direct consequence of the next result.

\begin{theorem}\label{t:main'}
The following assertions are valid:
\begin{enumerate}
\item[{$(i)$}]
Let $\mu<C_H+\lambda_1$ Then $\mathcal{N}=\{\sigma\leq \Lambda(p)\}$.
\item[{$(ii)$}]
Let $\mu=C_H+\lambda_1$. Then
$$\{\sigma<\Lambda(p)\}\cup\{\sigma=\Lambda(p),\,p\geq -1\}\subseteq\mathcal{N}\subseteq\{\sigma\le\Lambda(p)\}.$$
If $\,\Omega=S^{N-1}$ then $\mathcal{N}=\{\sigma<\Lambda(p)\}\cup\{\sigma=\Lambda (p),\,p\geq -1\}$.
\end{enumerate}
\end{theorem}

\noindent
Observe that in view of the scaling invariance of \eqref{e:main'} if
$w(x)$ is a solution to \eqref{e:main'} in $\C_\Omega^\rho$ then
$\tau^{\frac{2-s}{p-1}}w(\tau y)$ is a solution to \eqref{e:main'}
in $\C_\Omega^{\rho/\tau}$, for any $\tau>0$. So in what
follows, for $p\neq 1$, we confine ourselves to the study of
solutions to \eqref{e:main'} on $\C_\Omega^1$. For the same
reason, for $p\neq 1$ we may assume that $C=1$, when convenient.
In the remaining part of the paper we prove Theorem \ref{t:main'}.

\section{Improved Hardy inequality on cone--like domains}\label{s:Hardy}

Here and thereafter, for $0\le\rho<R\le+\infty$, we denote
$$\C_{\Omega}^{(\rho,R)}:=\{(r,\omega)\in\R^N:\:\omega\in\Omega,\:r\in(\rho,R)\},
\qquad\C_{\Omega}^\rho:=\C_{\Omega}^{(\rho,+\infty)},\qquad
\C_\Omega=\C_\Omega^0,$$ where $\Omega\subseteq S^{N-1}$ is a
subdomain of $S^{N-1}=\{x\in\R^N:|x|=1\}$.
%Accordingly, $\C_\Omega=\C_\Omega^0$ and $\C_{S^{N-1}}=\R^N\setminus\{0\}$.
We write $\Omega^{\prime}\Subset\Omega$ if $\Omega^{\prime}$ is a
smooth proper subdomain of $\Omega$ such that
$\Omega^{\prime}\neq\Omega$ and
$cl\,\Omega^{\prime}\subset\Omega$.
By $c,c_1,\dots$ we denote various positive constants exact values of which are irrelevant.
For two positive functions $\varphi_1$ and $\varphi_2$ we write $\varphi_1\asymp\varphi_2$
if there exist a constant $c\ge 1$ such that $c^{-1}\varphi_1\le\varphi_2\le c\varphi_1$.

Consider the linear equation
\begin{equation}\label{e:V}
-\Delta v-\frac{V(\omega)}{|x|^2}v=0\quad \mbox{in}\quad\C_\Omega^\rho,
\end{equation}
where $\Omega\subseteq S^{N-1}$ ($N\ge 2$) is a domain,
$V\in L^\infty(\Omega)$ and $\rho\ge 1$.
Recall that in the polar coordinates $(r,\omega)$ the operator $-\Delta-\frac{V(\omega)}{|x|^2}$ has the form
$$-\frac{\partial^2}{\partial r^2}-\frac{N-1}{r}\frac{\partial}{\partial r}+
\frac{1}{r^2}\left\{-\Delta_\omega-V(\omega)\right\},$$
where $\Delta_\omega$ denotes the Dirichlet Laplace--Beltrami operator on $\Omega$.
In what follows, $\lambda_{1,V}$ denotes the principal
%Dirichlet
eigenvalue
of the operator $-\Delta_\omega-V$ on $\Omega$.

The existence of positive solutions to \eqref{e:V} is equivalent
to the positivity of the quadratic form
$$\E_V(v,v):=\int_{\C_\Omega^\rho}|\nabla v|^2-\frac{V(\omega)}{|x|^2}v^2\,dx
%=\int_\rho^\infty\int_\Omega\left|\frac{\partial u}{\partial r}\right|^2+\frac{|\nabla_\omega u|^2}{r^2}
%-\frac{V(\omega)}{r^2}u^2\,d\omega\,dr,
\qquad(v\in H^1_c(\C_\Omega^\rho)\cap L^\infty_c(\C_\Omega^\rho)),$$
that corresponds to \eqref{e:V}, see \cite{Agmon} or Lemma \ref{l:A-ground} in Appendix.
Below we establish an improved Hardy--type inequality on cone--like domains,
which is appropriate for our purposes.
Similar inequalities were obtain recently on the ball and exterior domains
in \cite{Adimurthi,Barbatis}.

\begin{theorem}\label{l:Hardy}
%{\sf (Improved Hardy Inequality)}
The following inequality holds:
\begin{equation}\label{e:Hardy}
\E_V(v,v)\ge(C_H+\lambda_{V,1})\int_{\C_\Omega^\rho}\frac{v^2}{|x|^2}dx+
\frac{1}{4}\int_{\C_\Omega^\rho}\frac{v^2}{|x|^2\log^2|x|}dx,
\quad\forall\:v\in H^1_c(\C_\Omega^\rho)\cap L^\infty(\C_\Omega^\rho),
\end{equation}
where $C_H:=\left(\frac{N-2}{2}\right)^2$.
The constants $C_H+\lambda_{V,1}$ and $\frac{1}{4}$ are optimal in the sense that the inequality
\begin{equation}\label{e:Hardy-sharp}
\E_V(v,v)\ge\mu\int_{\C_\Omega^\rho}\frac{v^2}{|x|^2}dx+
\epsilon\int_{\C_\Omega^\rho}\frac{v^2}{|x|^2\log^2|x|}dx,
\quad\forall\:v\in H^1_c(\C_\Omega^\rho)\cap L^\infty(\C_\Omega^\rho),
\end{equation}
fails in any of the following two cases:
\begin{enumerate}
\item[$i)$]
$\mu=C_H+\lambda_{V,1}$ and $\,\epsilon>1/4$,
\item[$ii)$]
$\mu>C_H+\lambda_{V,1}$ and $\forall\,\epsilon\in\R$.
\end{enumerate}
%Note also that both constants $C_H+\lambda_1$ and $\frac{1}{4}$ do not depend on $\rho\ge 1$.
\end{theorem}

\begin{proof}%[Proof of Theorem \ref{l:Hardy}]
%Assume $\rho\ge 1$.
Let $\phi_\ast(r,\omega)=r^{\alpha_\ast}\log^{\frac{1}{2}}(r)\phi_{V,1}(\omega)$,
where $\phi_{1,V}>0$ is the eigenfunction of $-\Delta_\omega-V$,
that corresponds to $\lambda_{1,V}$. A direct computation shows that
$\phi_\ast\in H^1_{loc}(\C_\Omega^\rho)$ solves the equation
\begin{equation}\label{e:critical}
-\Delta v-\frac{V(\omega)}{|x|^2}v-\frac{C_H+\lambda_1}{|x|^2}v-\frac{1/4}{|x|^2\log^2|x|}v=0\quad\mbox{in }\:\C_\Omega^\rho.
\end{equation}
Thus the validity of \eqref{e:Hardy} follows from Lemma \ref{l:A-ground}.

Now we verify the optimality of constants in \eqref{e:Hardy}.
Fix $\rho\geq 1$ and let $R\gg\rho$.
Similarly to \cite{Adimurthi}, define a Lipschitz cut--off function
\begin{equation}\label{cut-off}
\theta_{\rho,R}(r):=\left\{
\begin{array}{cl}
0 & \text{for $r\le\rho$ and $r> R^2$},\\
r-\rho & \text{for $\rho<r\le \rho+1$},\\
1 & \text{for $\rho+1< r\leq R$}, \\
\frac{\log(R^2/r)}{\log(R)} & \text{for $R<r\leq R^2$}.\\
\end{array}
\right.
\end{equation}

$(i)$
We show that \eqref{e:Hardy-sharp} with $\mu=C_H+\lambda_{V,1}$ and $\epsilon>1/4$
fails on functions $\phi_\ast\theta_{\rho,R}\in H^1_c(\C_\Omega^\rho)\cap L^\infty(\C_\Omega^\rho)$.
By Lemma \ref{l:A-ground} direct computations give
\begin{eqnarray*}
&&
\int_{\C_\Omega^\rho}|\nabla(\phi_\ast\theta_{\rho,R})|^2-\frac{V(\omega)}{|x|^2}|\phi_\ast\theta_{\rho,R})|^2dx-
(C_H+\lambda_1)\int_{\C_\Omega^\rho}\frac{\phi_\ast^2\theta_{\rho,R}^2}{|x|^2}dx-
\epsilon\int_{\C_\Omega^\rho}\frac{\phi_\ast^2\theta_{\rho,R}^2}{|x|^2\log^2|x|}dx\\
&&\qquad\quad=
\int_\rho^\infty\int_\Omega\left|\nabla\left(\frac{\phi^2_\ast\theta^2_{\rho,R}}{\phi_\ast}\right)\right|^2
\phi_\ast^2\,d\omega\,r^{N-1}dr -\left(\epsilon-\frac{1}{4}\right)
\int_\rho^\infty\int_\Omega\frac{\phi^2_\ast\theta^2_{\rho,R}}{r^2\log^2(r)}\,d\omega\,r^{N-1}\,dr\\
%&&\qquad\quad=
%\int_\rho^\infty\left|\nabla\theta_{\rho,R}(r)\right|^2 r\log(r)\,dr
%\int_\Omega\phi_1^2d\omega\,
%-\left(\epsilon-\frac{1}{4}\right)\int_\rho^\infty\frac{\theta^2_{\rho,R}(r)}{r\log(r)}\,dr
%\int_\Omega\phi_1^2d\omega\,
%\\
&&\qquad\quad=
\int_\rho^{R^2}\left|\nabla\theta_{\rho,R}(r)\right|^2 r\log(r)\,dr
-\left(\epsilon-\frac{1}{4}\right)
\int_\rho^{R^2}\frac{\theta^2_{\rho,R}(r)}{r\log(r)}\,dr\\
&&\qquad\quad\leq
c_1+\int_R^{R^2}\frac{\log(r)}{r\log^2(R)}\,dr
-\left(\epsilon-\frac{1}{4}\right)\int_{\rho+1}^R \frac{1}{r\,\log r}\,dr\\
&&\qquad\quad\le
c_2-\left(\epsilon-\frac{1}{4}\right)\log\log(R)\to-\infty
\quad\mbox{as }\:R\to\infty.
\end{eqnarray*}
Observe that the result does not depend on the initial choice of $\rho\ge 1$.

$(ii)$
Choosing $\phi_\ast(\rho,\omega)=r^{\alpha_\ast}\phi_{V,1}(\omega)$,
one can verify that \eqref{e:Hardy-sharp} with $\mu>C_H+\lambda_{V,1}$ and any $\epsilon\in\R$ fails on functions
$\theta_{\rho,R}\phi_\ast\in H^1_c(\C_\Omega^\rho)\cap L^\infty(\C_\Omega^\rho)$ for large $R\gg\rho$.
We omit the details.
\end{proof}

Optimality of Improved Hardy Inequality \eqref{e:Hardy}, via Corollary \ref{APA}, implies
the following nonexistence result, which is one of the main tools
in our proofs of nonexistence of positive solutions to semilinear equation \eqref{e:main'}.

\begin{corollary}\label{l:Nonexist}
%{\sc (Nonexistence Lemma)}
Equation \eqref{e:V} has a positive super-solution
if and only if $C_H+\lambda_{V,1}\ge 0$.
\end{corollary}

\begin{remark}\label{r:Nonexist}
In particular, if $V(\omega)\equiv\mu$ then equation \eqref{e:V}
has a positive super-solution if and only if $\mu\le C_H+\lambda_1$.
\end{remark}

\section{Asymptotics of positive super-solutions to $-\Delta v-\frac{V(\omega)}{|x|^2}v=0$}
\label{s:asymptotic}

According to Lemma \ref{l:Nonexist},
equation \eqref{e:V} admits positive super-solutions if and only if $C_H+\lambda_{V,1}\ge 0$.
In this section, by constructing appropriate comparison functions,
we obtain sharp two-sided bounds on the growth at infinity of super-solutions to \eqref{e:V}.

Throughout this section $(\lambda_{V,k})_{k\in\N}$ denotes the sequence of Dirichlet eigenvalues of the operator
$-\Delta_\omega-V$ in $L^2(\Omega)$,
$$\lambda_{V,1}<\lambda_{V,2}\le\dots\le\lambda_{V,k}\le\dots.$$
By $(\phi_{V,k})_{k\in\N}$ we denote the corresponding orthonormal basis
of eigenfunctions in $L^2(\Omega)$, with the positive principal eigenfunction $\phi_{V,1}>0$.
%Note that by the standard elliptic regularity $\phi_{V,k}\in L^\infty(\Omega)$.
If $V=0$ and there is no ambiguity we simply write $\lambda_k$ and $\phi_k$.

Let $\psi\in L^2(\Omega)$. Then
\begin{equation}\label{Fourier}
\psi=\sum^\infty_{k=1}\psi_k\phi_{V,k},\qquad
\mbox{where}\quad\psi_k=\int_\Omega\psi\phi_{V,k}\,d\omega,
\end{equation}
and the series converges in $L^2(\Omega)$ with $\|\psi\|^2_{L^2}=\sum_{k=1}^\infty\psi_{k}^2$.
If, in addition, $\psi\in H^1_0(\Omega)$ then \eqref{Fourier} converges in $H^1_0(\Omega)$
with $\|\nabla\psi\|^2_{L^2}\asymp\sum_{k=1}^\infty\lambda_{V,k}\psi_{k}^2$.
%In particular, if $V=\mu$ where $\mu$ is a constant, then
%$\lambda_{V,k}=\lambda_k-\mu$, where $(\lambda_k)_{k\in\N}$ is a
%sequence of Dirichlet eigenvalues of the Laplace-Beltrami operator $-\Delta_\omega$ on $\Omega$.
In what follows we use the following simple observation.

\begin{lemma}\label{l-series}
Let $\psi\in C^\infty_c(\Omega)$. Then the Fourier series \eqref{Fourier} converges in $L^\infty(\Omega)$.
\end{lemma}

\proof
Observe that
$\|\phi_k\|_{L^\infty}\le c\lambda_{V,k}^{\frac{N-1}{4}}$,
by the standard elliptic estimates of eigenfunctions
of the Dirichlet Laplace--Beltrami operator $-\Delta_\omega-V$ on $\Omega\subseteq S^{N-1}$,
see, e.g. \cite[p.172]{Chavel}.
Choose $b>(N-1)/2$ and $a=(N-1)/4-b$.  Then
$$\sum^\infty_{k=1}|\psi_k||\phi_{V,k}|\le c\sum^\infty_{k=1}|\psi_k|\lambda_{V,k}^{\frac{N-1}{4}}
\le c\left(\sum^\infty_{k=1}\lambda_{V,k}^{2a}\right)^{1/2}
\left(\sum^\infty_{k=1}|\psi_k|^2\lambda_{V,k}^{2b}\right)^{1/2}
<\infty.$$
Here the first series converges due to the classical spectral asymptotics
$\lambda_{V,k}\asymp k^{\frac{2}{N-1}}$.
The second series series converges by a standard spectral argument,
taking into account that $\psi\in\cap_{m\in\N}D((-\Delta_\omega)^m)$,
where $D((-\Delta_\omega)^m)$ is the domain of the $m$--the power of the Dirichlet Laplace--Beltrami
operator $-\Delta_\omega$ on $\Omega$.
\qed

If $C_H+\lambda_1\ge 0$ then the roots of the quadratic equation
\begin{equation}\label{e:alpha-k}
\alpha(\alpha+N-2)=\lambda_{V,k}
\end{equation}
are real, for each $k\in\N$. Denote these roots by $\alpha_{V,k}^-\le\alpha_{V,k}^+$.
%$$\alpha_k^{\pm}:=-\frac{N-2}{2}\pm\sqrt{\frac{(N-2)^2}{4}+(\lambda_k-\mu)}.$$
If $C_H+\lambda_{V,1}=0$ and $k=1$ then \eqref{e:alpha-k} has the unique root,
denoted by $\alpha_\ast:=\alpha_{V,1}^\pm=\frac{2-N}{2}$.

For a positive function $u\in H^1_{loc}(\C_\Omega^1)$ and a
subdomain $\Omega^{\prime}\subseteq\Omega$, denote
$$m_{u}(R,\Omega^{\prime}):=\inf_{\C_{\Omega^{\prime}}^{(R/2,R)}}u,\qquad
M_{u}(R,\Omega^{\prime}):=\sup_{\C_{\Omega^{\prime}}^{(R/2,R)}}u.$$

Our main result in this section reads as follows.
\begin{theorem}\label{main:bound}
Let $u\in H^1_{loc}(\C_\Omega^1)$ be a positive super-solution to \eqref{e:V}.
Then for every proper subdomain $\Omega^{\prime}\Subset\Omega$ and $R\gg 1$ the following hold:
\begin{enumerate}
\item[$(i)$]
if $C_H+\lambda_{V,1}>0$ then
\begin{equation}\label{e:main:bound}
c_1R^{\alpha_{V,1}^-}\le m_u(R,\Omega^{\prime})\le c_2 R^{\alpha_{V,1}^+},
\end{equation}
\item[$(ii)$]
if $C_H+\lambda_{V,1}=0$ then
\begin{equation}\label{e:main:bound-eps}
c_1 R^{\alpha_\ast}\le m_u(R,\Omega^{\prime})\le c_2 R^{\alpha_\ast}\log(R).
\end{equation}
\end{enumerate}
\end{theorem}

\begin{remark}
The above estimates are sharp, as one sees comparing with the explicit solutions
$r^{\alpha_{V,1}^\pm}\phi_{V,1}$ in the case $(i)$ and
$r^{\alpha_\ast}\phi_{V,1}$ and $r^{\alpha_\ast}\log(r)\phi_{V,1}$ in the case $(ii)$.
\end{remark}

\begin{remark}\label{r:scaling}
Equation \eqref{e:V} is invariant with respect to scaling.
Namely, if $v(x)$ is a (super)\,solution to \eqref{e:V} in $\C_\Omega^\rho$,
then $v(\tau x)$ is a solution to \eqref{e:V} in $\C_\Omega^{\tau\rho}$, for any $\tau>0$.
Therefore, in what follows we may consider \eqref{e:V} in $\C^\rho_\Omega$
with a conveniently fixed radius $\rho\ge 1$.
\end{remark}

\begin{remark}\label{r:Harnack}
The scaling invariance implies that positive (super)\,solutions to \eqref{e:V}
satisfy the Harnack inequalities with $r$--independent constants.
More precisely, if $u>0$ is a super-solution to \eqref{e:V}
then the weak Harnack inequality reads as
\begin{equation}\label{h:weak}
\int_{\C_{\Omega^{\prime}}^{(R/2,R)}}u\,dx\le C_w\,R^N
m_u(R,\Omega^{\prime}),
\end{equation}
where $C_w=C_w(\Omega^\prime)>0$ does not depend on $R\gg 1$.
Similarly, if $u>0$ is a solution to \eqref{e:V} then
by the strong Harnack inequality
\begin{equation}\label{h:strong}
M_u(R,\Omega^{\prime})\le C_s\,m_u(R,\Omega^{\prime}),
\end{equation}
where $C_s=C_s(\Omega^\prime)>0$ is independent of $R\gg 1$.
\end{remark}

In the remaining part of the section we prove Theorem \ref{main:bound}.
Our proof relies on the Comparison Principle in the extended Dirichlet spaces
associated to \eqref{e:V} (see Appendix \ref{A1}).
The cases $(i)$ and $(ii)$ are considered separately.

\subsection{Case $C_H+\lambda_{V,1}>0$}

The above condition is assumed throughout this subsection.
In this case Hardy Inequality \eqref{e:Hardy} implies that the form
$\E_V$ satisfies $\lambda$--property \eqref{e:lambda} with
$\lambda(x)=\frac{C_H+\lambda_{V,1}}{|x|^2}$.
Hence the extended Dirichlet space $\D(\E_V,\C_\Omega^2)$ is well--defined
(see Appendix \ref{A1}) and the Comparison Principle (Lemma \ref{l:WCP}) is valid.
Moreover, \eqref{e:Hardy} implies that
$$\D(\E_V,\C_\Omega^2)=D^1_0(\C_\Omega^2),$$
where $D^1_0(\C_\Omega^2)$ is the usual homogeneous Sobolev space,
defined as the completion of $C^\infty_c(\C_\Omega^2)$
with respect to the Dirichlet norm $\|\nabla u\|_{L^2}$.

\paragraph{Lower estimate.}
Fix a smooth proper subdomain $\Omega^{\prime}\Subset\Omega$ and
a function $0\lneq\psi\in C^\infty_c(\Omega^{\prime})$.
%Without loss of generality we may assume that $\psi<\phi_1$.
For $(r,\omega)\in\C_\Omega^2$ and $k\in\N$ set
\begin{equation}\label{v-k}
v_k(r,\omega):=\eta_k(r)\phi_{V,k}(\omega),
\qquad\mbox{where}\quad\eta_k(r):=\left(\frac{r}{2}\right)^{\alpha_{V,k}^-}.
\end{equation}
Define the \textit{comparison function} $v_\psi$ by
\begin{equation}\label{e:v-psi}
v_\psi:=\sum^\infty_{k=1}\psi_k v_k,
\end{equation}
where $\psi_k$ are the Fourier coefficients of $\psi$ as in \eqref{Fourier}.
Thus $v_\psi(2,\omega)=\psi(\omega)$.
A direct computation verifies that $v_\psi\in H^1_{loc}(\C_\Omega^2)$
is a solution to \eqref{e:V} in $\C_\Omega^2$.

\begin{lemma}\label{l:Minimal}
Let $0<u\in H^1_{loc}(\C_\Omega^1)$ be a super-solution to \eqref{e:V} in $\C_\Omega^1$.
Then
$$u\geq cv_\psi\quad in\quad \C^2_{\Omega^{\prime}}.$$
\end{lemma}

\proof
Choose a function $\theta(r)\in C^{0,1}[2,+\infty)$ such that
$0\leq\theta(r)\leq 1$, $\theta(2)=1$ and $\theta(r)=0$ for $r\ge 3$.
Set $\tilde v_k:=v_k-\theta\phi_{V,k}$.
By direct computations,
\begin{equation}\label{e:v-k}
\E_V(\tilde v_k,\tilde v_k)\le c_1|\alpha_{V,k}^-|+c_2\qquad\mbox{and}\qquad
\E_V(\tilde v_k,\tilde v_l)=0\quad\mbox{for }\:l\neq k.
\end{equation}
Then it is straightforward that $\tilde v_k\in D^1_0(\C_\Omega^2)$.
Consider the function
$$\tilde v_\psi:=\sum^\infty_{k=1}\psi_k\tilde v_k.$$
By \eqref{e:v-k} and taking into account that $|\alpha_k^-|\asymp\sqrt{\lambda_k}$ we obtain
\begin{eqnarray*}
\E_V(\tilde v_\psi,\tilde v_\psi)&\le&
\sum^\infty_{k=1}\psi_k^2(c_1|\alpha_k^-|+c_2)\le
c_3\left(\sum^\infty_{k=1}\psi_k^2\lambda_k\right)^{1/2}\left(\sum^\infty_{k=1}\psi_k^2\right)^{1/2}+
c_2 \sum^\infty_{k=1}\psi_k^2\\
&=&
c_3\|\nabla_\omega\psi\|_{L^2}\|\psi\|_{L^2}+c_2\|\psi\|^2_{L^2}.
\end{eqnarray*}
Hence $\tilde v_\psi=v_\psi-\theta\psi\in D^1_0(\C_\Omega^2)$.

Now observe that by the weak Harnack inequality \eqref{h:weak}
there exists $\delta>0$ such that
$$u>\delta\quad\mbox{in}\quad\C^{(2,3)}_{\Omega^{\prime}}.$$
Fix $c>0$ such that $c\psi<\delta$ in $\Omega^{\prime}$. Thus
$u>c\theta\psi$ in $\C_{\Omega}^2$. Represent
$$u-cv_\psi=(u-c\psi\theta)-c\tilde v_\psi,$$
where $\tilde v_\psi\in D^1_0(\C_\Omega^2)$
and notice that $u-v_\psi$ is a super-solution to \eqref{e:V}.
By Lemma \ref{l:WCP} we conclude that $u-c\psi\theta\geq c\tilde v_\psi$,
that is $u\geq c v_\psi$ in $\C_\Omega^2$.
\qed

\begin{lemma}\label{l:v-psi}
$m_{v_{\psi}}(R,\Omega^{\prime})\asymp R^{\alpha_1^-}$ as $R\to\infty$.
\end{lemma}

\begin{proof}
Choosing $u=r^{\alpha_{V,1}^-}\phi_1$ as a (super)\,solution in Lemma \ref{l:Minimal}
we immediately conclude that
$$m_{v_{\psi}}(R,\Omega^{\prime})\le cR^{\alpha_{V,1}^-}\quad\mbox{for }\:R\gg 2.$$
To obtain the reverse inequality, note that $v_\psi$ as $v_\psi=\psi_1 v_1+w_\psi$.
%$$v_\psi=\psi_1 v_1+w_\psi,\qquad\mbox{where}\quad
%w_\psi=\sum_{k=2}^\infty\psi_k \phi_k(\omega)\left(\frac{r}{2}\right)^{\alpha_{V,k}^-}.$$
Then by Lemma \ref{l-series} we obtain the uniform bound
\begin{equation}\label{e:series}
|w_\psi(r,\omega)|\le
\eta_2(r)\sum_{k=2}^\infty|\psi_k||\phi_k(\omega)|\le c r^{\alpha_{V,2}^-}.
\end{equation}
Note that $\phi_{V,1}>\delta$ in $\Omega^{\prime}$, for some $\delta>0$.
We conclude that
$$m_{v_\psi}(\Omega^{\prime},R)\ge c_2 R^{\alpha_{V,1}^-}-c_3 R^{\alpha_{V,2}^-}\quad\mbox{for }\:R\gg 4.$$
This completes the proof, since $\alpha_{V,2}^-<\alpha_{V,1}^-<0$.
\end{proof}

Combining Lemmas \ref{l:Minimal} and \ref{l:v-psi} we obtain the lower bound in \eqref{e:main:bound}.

\paragraph{Upper estimate.}

Fix a subdomain $\Omega^{\prime}\subseteq\Omega$ and a function
$0\leq\psi\in C^\infty_c(\Omega^{\prime})$. Let $R\ge 4$. For
$(r,\omega)\in\C_\Omega^{(1,R)}$ and $k\in\N$ define
\begin{equation}\label{v-kR}
v_{k,R}(r,\omega):=\eta_{k,R}(r)\phi_{V,k}(\omega),
\qquad\mbox{where}\quad
\eta_{k,R}(r):=\left\{\frac{r^{\alpha_{V,k}^+}-r^{\alpha_{V,k}^-}}
{R^{\alpha_{V,k}^+}-R^{\alpha_{V,k}^-}}\right\}.
\end{equation}
Let $\theta:[0,1]\to\R$ be a smooth function such that
$0\le\theta\le 1$, $\theta(1)=1$ and $\theta(\xi)=0$ for $\xi\in[0,1/2]$.
For $r\in[R/2,R]$ set $\theta_R(r):=\theta(r/R)$.
%Thus $\theta_R\psi\in \D(\E_\mu,\C_{\Omega^{\prime}}^{(R/2,R)})$.
A direct computation verifies that $v_{k,R}$ is a solution to the problem
\begin{equation}\label{e-2R}
\left(-\Delta -\frac{V(\omega)}{|x|^2}\right)v=0,\qquad v-\theta_R\phi_{V,k}\in H^1_0(\C_\Omega^{(1,R)}).
\end{equation}
Let
\begin{equation*}\label{v-psiR}
v_{\psi,R}:=\sum^\infty_{k=1}\psi_k\,v_{k,R}
\end{equation*}
where $\psi_k$ are the Fourier coefficients of $\psi$ as in \eqref{Fourier}.
Thus $v_{\psi,R}\in H^1_{loc}(\C_\Omega^{1,R})$ is a solution to (\ref{e-2R})
and $v_{\psi,R}(R,\omega)=\psi(\omega)$.

Fix a compact $K_0\subset\C^{(2,3)}_\Omega$.
Define the \textit{comparison function} $v_{\psi,R}$ by
$$\tilde v_{\psi,R}:=\frac{v_{\psi,R}}{\inf_{K_0}v_{\psi,R}}.$$
Then $\inf_{K_0}\tilde v_{\psi,R}= 1$.
%and $\tilde v_{\psi,R}$ solves \eqref{e-2R}.
Note that the construction of $\tilde v_{\psi,R}$ depends only on the choice of $K_0$, $\psi$ and $R$.
The following lemma is a weak version of the Phragmen--Lindel\"of comparison principle
adopted to our framework.

\begin{lemma}\label{l:Phragmen}
%{\sf (Phragmen--Lindel\"of type comparison principle)}
Let $0<u\in H^1_{loc}(\C_\Omega^1)$ be a super-solution to \eqref{e:V} in $\C_\Omega^1$.
Then
$$m_{u}(R,\Omega^{\prime})\le cM_{\tilde v_{\psi,R}}(R,\Omega^{\prime}),\qquad R\ge 4.$$
\end{lemma}

\proof
Set $\delta_R:=\inf_{K_0}v_{\psi,R}$. For a contradiction assume
that for any $c>0$ there exists $R\ge 4$ such that
$$u\ge c\tilde v_{\psi,R}=\frac{c}{\delta_R}v_{\psi,R}
\quad\text{in }\:\C_{\Omega^{\prime}}^{(R/2,R)}.$$ Let $\psi_R>0$
be the unique solution to the problem
$$-\Delta v-\frac{V(\omega)}{|x|^2}v = 0,\qquad v-\theta_R\psi\in H^1_0(\C_{\Omega^{\prime}}^{(R/2,R)}).$$
Then clearly
$$\left(-\Delta -\frac{V(\omega)}{|x|^2}\right)(v_{\psi,R}-\psi_R)=0\quad\text{in }\:\C_{\Omega^{\prime}}^{(R/2,R)}.$$
Observe that $v_{\psi,R}>0$ in
$\C_{\Omega}^{(1,R)}\setminus\C_{\Omega^{\prime}}^{(R/2,R)}$.
Hence $(v_{\psi,R}-\psi_R)^-\in D^1_0(\C_{\Omega^{\prime}}^{R/2,R})$.
By Lemma \ref{l:WMP} we conclude that
$$v_{\psi,R}\ge\psi_R\quad\text{in }\:\C_{\Omega^{\prime}}^{(R/2,R)}.$$
Let $\bar\psi_R$ denote the function $\psi_R$, extended to $\C_\Omega^{(1,R)}$ by zero.
Therefore
$$\left(-\Delta-\frac{V(\omega)}{|x|^2}\right)(u-c \tilde v_{\psi,R})=
\left(-\Delta-\frac{V(\omega)}{|x|^2}\right)\left((u-\frac{c}{\delta_R}\bar\psi_R)-\frac{c}{\delta_R}
(v_{\psi,R}-\bar\psi_R)\right)\ge 0
\quad\mbox{in }\:\C_\Omega^{(1,R)}.$$ Then Lemma \ref{l:WCP} implies that
\begin{equation*}
u\ge c\tilde v_{\psi,R}\quad\mbox{in }\:\C_\Omega^{(1,R)}.
\end{equation*}
Since $c>0$ is arbitrary, we conclude that $\inf_{K_0}u=+\infty$.
Hence, by weak Harnack inequality \eqref{h:weak}, $u\equiv+\infty$ in $\C_\Omega^1$,
which is a contradiction.
\qed

\begin{lemma}\label{l:w-psiR}
$M_{\tilde v_{\psi,R}}(R,\Omega^{\prime})\asymp R^{\alpha_{V,1}^+}$ as $R\to\infty$.
\end{lemma}

\begin{proof}
Choosing $u:=r^{\alpha_{V,1}^+}\phi_1$ as a (super)\,solution in Lemma \ref{l:Phragmen} we conclude that
$$M_{\tilde v_{\psi,R}}(R,\Omega^{\prime})\ge cR^{\alpha_{V,1}^+},\qquad R\gg 1.$$
Now we estimate $M_{\tilde v_{\psi,R}}(R,\Omega^{\prime})$ from above.

First, observe that Lemma \ref{l:WCP}, Lemma \ref{l-series}
and the arguments, similar to those in Lemma \ref{l:Minimal} imply the upper bound
$$v_{\psi,R}(r,\omega)\le c\eta_{1,R}(r)\phi_{V,1}(\omega)\quad\mbox{in }\:\C_\Omega^{(1,R)},$$
where $c>0$ is chosen so that $\psi\le c\phi_{V,1}$ in $\Omega$.
Clearly, if $\alpha_{V,1}^+\ge 0$ then $\eta_{1,R}(r)\le 1$.
However, if $\alpha_{V,1}^+<0$ then $\eta_{1,R}(r)$ attains its maximum at $r_\ast\in(1,R)$
%$$r_\ast=\left(\frac{\alpha_k^+}{\alpha_k^-}\right)^{\frac{1}{\alpha_k^+-\alpha_k^-}}>1,$$
with $\eta_{1,R}(r_\ast)\to\infty$ as $R\to\infty$.
Nevertheless, one can readily verify that%it is easily seen that
$$\max_{r\in[R/2,R]}\eta_{k,R}(r)\le\max\{1,2^{-\alpha_{V,1}^+}\},\qquad R\gg 1.$$
Therefore
$$M_{v_{\psi,R}}(\Omega,R)\le c_1,\qquad R\gg 1.$$
To estimate $\inf_{K_0}v_{\psi,R}$ from below, note that %$v_{\psi,R}$ as
$$v_{\psi,R}=\psi_1 v_{1,R}+w_{\psi,R}\,,
\qquad\mbox{where}\quad w_{\psi,R}=\sum_{k=2}^\infty\psi_k v_{k,R}\,.$$
Then by Lemma \ref{l-series} similarly to \eqref{e:series} we obtain
$$\sup_{K_0}|w_{\psi,R}|\le
\max_{r\in(2,3)}\eta_k(r)
\sum_{k=2}^\infty |\psi_k||\phi_{V,k}(\omega)|\le
\frac{c_1}{R^{\alpha_{V,2}^+}-R^{\alpha_{V,2}^-}}.$$
We conclude that
$$\inf_{K_0}v_{\psi,R}\ge
\inf_{K_0}\psi_1 v_{1,R}(r)
-\sup_{K_0}|w_{\psi,R}|\ge\frac{c_2}{R^{\alpha_1^+}}-\frac{c_3}{R^{\alpha_2^+}}.$$
This completes the proof since $\alpha_{V,2}^+>\alpha_{V,1}^+$.
\end{proof}

Combining Lemmas \ref{l:Phragmen} and \ref{l:w-psiR} we obtain the upper bound in \eqref{e:main:bound}.

%\begin{remark}
%Note that the choice $\Omega^{\prime}=\Omega$ is allowed in Lemma \ref{l:Phragmen} and \ref{l:w-psiR}.
%\end{remark}

\subsection{Case $C_H+\lambda_{V,1}=0$}

The above condition is assumed throughout this subsection.
%\paragraph{The space $\tilde D^1_0(\C_\Omega^\rho)$.}
Let $\rho\ge 1$.
Then Hardy's Inequality \eqref{e:Hardy} implies that the form $\E_V$
satisfies the $\lambda$--property \eqref{e:lambda} with $\lambda(x)=\frac{1/4}{|x|^2\log^2|x|}$.
Hence the extended Dirichlet space $\D(\E_V,\C_\Omega^\rho)$ is well--defined
(see Appendix \ref{A1}), and in particular, the Comparison Principle (Lemma \ref{l:WCP}) is valid.
We denote
$$\tilde D^1_0(\C_\Omega^\rho):=\D(\E_V,\C_\Omega^\rho).$$
%Then Hardy's Inequality \eqref{e:Hardy} implies that
%$$\D(\E_{C_H+\lambda_1,\epsilon},\C_\Omega^\rho)=\tilde D^1_0(\C_\Omega^\rho),\qquad
%\forall\:\epsilon\in(0,1/4).$$
The space $\tilde D^1_0(\C_\Omega^2)$ is larger then $D^1_0(\C_\Omega^2)$ (cf.\ \cite{Barbatis,VZ}).
In order to see this, for $\beta\in[0,1]$ consider
\begin{equation}\label{v-beta}
v_\beta(r,\omega):=r^{\alpha_\ast}\log^\beta(r)\phi_{V,1}(\omega).
\end{equation}
Clearly, $v_{\beta}\in C^\infty_{loc}(\C_\Omega^\rho)$
but $\nabla v_\beta\not\in L^2(\C_\Omega^\rho)$.
%A direct computation verifies that $v_\beta$ solves \eqref{e:V}.
Let $\theta(r)\in C^{0,1}[\rho,+\infty)$ be such that
$0\leq\theta(r)\leq 1$, $\theta(\rho)=1$ and $\theta(r)=0$ for $r\ge\rho+1$.
\begin{lemma}\label{p:A1}
$\:v_\beta-\theta\phi_{V,1}\in\tilde D^1_0(\C_\Omega^\rho)$ for each $\beta\in[0,1/2)$.
\end{lemma}
\begin{proof}
Define the cut--off function $\theta_R(r)\in C_{c}^{\,0,1}(\C_\Omega^1)$ by
$$\theta_R(r):=\left\{ \begin{array}{ll}
1, & 1\leq r\leq R, \\
\frac{\log(R^2/r)}{\log R}, & R\leq r\leq R^2,\\
0, & r\geq R^2.
\end{array}
\right.$$
Let $w_R:=\theta_R(v_\beta-\theta\phi_{V,1})$.
According to Lemma \ref{l:A-ground}, one can represent $\E_V(v_R)$ as
\begin{eqnarray}\label{e:A1}
\E_V(w_R,w_R)
&=&\int_\rho^\infty\int_\Omega\left|\nabla\left(\frac{w_R}{v_0}\right)\right|^2 v_0^2\,d\omega\,r^{N-1}dr
%=\int_\rho^\infty\left|\nabla\left(\log^\beta(r)\theta_R(r)\right)\right|^2 r \,dr\,
%\int_\Omega\phi_1^2\,d\omega
\nonumber\\
&=&\int_\rho^\infty\left|\nabla\left(\log^\beta(r)\theta_R(r)\right)\right|^2 r \,dr
\le c_1+c_2\log^{2\beta-1}(R)\le c.\nonumber
\end{eqnarray}
Hence $\E_V(w_{R_n},w_{R_n})$ is a Cauchy sequence, for an appropriate choice of $R_n\to\infty$.
Since $(w_{R_n})\subset C_{loc}^{0,1}(\C_\Omega^\rho)$ converges pointwise to the function $v_\beta$,
the assertion follows.
\end{proof}

Now we are in a position to prove \eqref{e:main:bound-eps}.

\paragraph{Lower estimate.}
As before, fix a proper smooth subdomain
$\Omega^{\prime}\Subset\Omega$ and a function $0\lneq\psi\in
C^\infty_c(\Omega^{\prime})$.
%Without loss of generality we may assume that $\psi<\phi_1$.
For $(r,\omega)\in\C_\Omega^2$ and $k\in\N$ set
$$v_\ast(r,\omega):=c_\ast r^{\alpha_\ast}\phi_{V,1}(\omega),$$
where $c_\ast>0$ chosen so that $v_\ast(2,\omega)=\phi_{V,1}(\omega)$.
Clearly $v_\ast\in H^1_{loc}(\C_\Omega^2)$ is a solution to \eqref{e:V} in $\C_\Omega^2$.

Define the \textit{comparison function} $v_\psi$ by
\begin{equation}\label{e:v-psi-eps}
v_\psi:=\psi_1 v_\ast +\sum^\infty_{k=2}\psi_k v_k,
\end{equation}
where $\psi_k$ are the Fourier coefficients of $\psi$ as in
\eqref{Fourier} and $v_k$ with $k\ge 2$ are defined by
\eqref{v-k}. Thus $v_\psi(2,\omega)=\psi(\omega)$. Observe that
for $k\ge 2$ the functions $v_k$ are solutions to \eqref{e:V} in $\C_\Omega^2$.
Hence $v_\psi\in H^1_{loc}(\C_\Omega^2)$ is a solution to \eqref{e:V} in $\C_\Omega^2$.

\begin{lemma}\label{l:Minimal-eps}
Let $0<u\in H^1_{loc}(\C_\Omega^1)$ be a super-solution to \eqref{e:V} in $\C_\Omega^1$.
Then
$$u\geq cv_\psi\quad in\quad \C^2_{\Omega^{\prime}}.$$
\end{lemma}

\begin{proof}
Similar to the proof of Lemma \ref{l:Minimal}.
\end{proof}

\begin{lemma}\label{l:v-psi-eps}
$m_{v_{\psi}}(R,\Omega^{\prime})\asymp R^{\alpha_\ast}$ as $R\to\infty$.
\end{lemma}

\begin{proof}
Similar to the proof of Lemma \ref{l:v-psi}.
\end{proof}

Combining Lemmas \ref{l:Minimal-eps} and \ref{l:v-psi-eps}
we obtain the lower bound in \eqref{e:main:bound-eps}.

\paragraph{Upper estimate.}
Fix a subdomain $\Omega^{\prime}\subseteq\Omega$ and a function
$0\leq\psi\in C^\infty_c(\Omega^{\prime})$. Let $R\ge 4$. For
$(r,\omega)\in\C_\Omega^{(1,R)}$ and $k\in\N$ define
\begin{equation}\label{v-kR-eps}
v_{\ast,R}(r,\omega):=\eta_{\ast,R}(r)\phi_{V,1}(\omega),
\quad\mbox{where}\quad
\eta_{\ast,R}(r):=\frac{\log(r)}{\log(R)}\left(\frac{r}{R}\right)^{\alpha_\ast}.
\end{equation}
Let $\theta:[0,1]\to\R$ be a smooth function such that
$0\le\theta\le 1$, $\theta(1)=1$ and $\theta(\xi)=0$ for $\xi\in[0,1/2]$.
For $r\in[R/2,R]$ set $\theta_R(r):=\theta(r/R)$.
%Thus $\theta_R\psi\in \D(\E_\mu,\C_{\Omega^{\prime}}^{(R/2,R)})$.
A direct computation verifies that $v_{\ast,R}$ and $v_{k,R}$
($k\ge 2$), defined by \eqref{v-k} are solutions to the problems
\begin{equation}\label{e-2R-eps}
\left(-\Delta -\frac{V(\omega)}{|x|^2}\right)v=0,
\qquad v-\theta_R\phi_{V,k}\in H^1_0(\C_\Omega^{(1,R)}).
\end{equation}
Let
\begin{equation*}\label{v-psiR-eps}
v_{\psi,R}:=\psi_1\,v_{\ast,R}+\sum^\infty_{k=2}\psi_k\,v_{k,R}
\end{equation*}
where $\psi_k$ are the Fourier coefficients of $\psi$ as in
\eqref{Fourier}. Thus $v_{\psi,R}\in H^1_{loc}(\C_\Omega^{1,R})$
is a solution to (\ref{e-2R-eps}) and
$v_\psi(R,\omega)=\psi(\omega)$.

Fix a compact $K_0\subset\C^{(2,3)}_\Omega$.
Define the \textit{comparison function} $\tilde v_{\psi,R}$ by
$$\tilde v_{\psi,R}:=\frac{v_{\psi,R}}{\inf_{K_0}v_{\psi,R}}.$$
%so $\inf_{K_0}\tilde v_{\psi,R}=1$.

\begin{lemma}\label{l:Phragmen-eps}
%{\sf (Phragmen--Lindel\"of type comparison principle)}
Let $0<u\in H^1_{loc}(\C_\Omega^1)$ be a super-solution to \eqref{e:V} in $\C_\Omega^1$.
Then
$$m_{u}(R,\Omega^{\prime})<c M_{\tilde v_{\psi,R}}(R,\Omega^{\prime}),\qquad R\ge 4.$$
\end{lemma}

\begin{proof}
Similar to the proof of Lemma \ref{l:Phragmen}.
\end{proof}

\begin{lemma}\label{l:w-psiR-eps}
$M_{\tilde v_{\psi,R}}(R,\Omega^{\prime})\asymp R^{\alpha_\ast}\log(R)$ as $R\to\infty$.
\end{lemma}

\begin{proof}
Similar to the proof of Lemma \ref{l:w-psiR}.
\end{proof}

Combining Lemmas \ref{l:Phragmen-eps} and \ref{l:w-psiR-eps}
we obtain the upper bound in \eqref{e:main:bound-eps}.

\subsection{Auxiliary linear equation}

In this subsection we consider the inhomogeneous linear equation
\begin{equation}\label{e-log-sigma}
-\Delta
w-\frac{C_H+\lambda_1}{|x|^2}w=\frac{\psi(\omega)}{|x|^{2-\alpha_\ast}\log^\sigma|x|}
\quad\mbox{in}\quad\C_\Omega^\rho,
\end{equation}
where $\alpha_\ast=\frac{2-N}{2}$, $\sigma>0$, $\rho>\exp(1)$ and $0\lneq \psi\in C^\infty_c(\Omega)$.
\begin{lemma}\label{l-log-sigma}
Equation \eqref{e-log-sigma} has no positive super-solution for $\sigma\le 1$.
\end{lemma}
\begin{proof}
Without loss of generality we may assume $\sigma=1$.
For each $R\gg\rho$ we are going to construct a barrier $w_{\psi,R}>0$ that solves the problem
\begin{equation}\label{e-loglog-sum}
-\Delta w-\frac{C_H+\lambda_1}{|x|^2}w=\frac{\psi(\omega)}{|x|^{2-\alpha_\ast}\log^\sigma|x|},
\qquad w\in H^1_0(\C_\Omega^{(\rho,R)}),
\end{equation}
and blows up on a fixed compact $K\subset \C_\Omega^{\rho}$ as $R\to\infty$.
Then by Lemma \ref{l:WCP},
$$u\ge w_{\psi,R}\qquad\mbox{in }\C_\Omega^{(\rho,R)}.$$
Therefore we conclude that $u\equiv +\infty$ in $K$, which is a contradiction.

To construct such $w_{\psi,R}$, consider the boundary value problem
\begin{equation}\label{ode-eta}
-\eta_k^{\prime\prime}-\frac{N-1}{r}\eta_k^\prime-\frac{C_H}{r^2}\eta_k+\frac{\delta_k^2}{r^2}\eta_k=
\frac{1}{r^{2-\alpha_\ast}\log^\sigma(r)},\qquad\eta(\rho)=\eta(R)=0,
\end{equation}
where $\delta_k:=\sqrt{\lambda_k-\lambda_1}$, and $k\in\N$.
For $k=1$, the solution to \eqref{ode-eta} is given by
$$\eta_{1,R}(r)=r^{\alpha_\ast}
\left(A_{1,R}+B_{1,R}\log(r)+\log(r)\log\log(r)\right),$$
where
$$A_{1,R}=\frac{\log(R)\log(\rho)(\log\log(\rho)-\log\log(R))}{\log(R)-\log(\rho)},\quad
B_{1,R}=\frac{\log(R)\log\log(R)-\log(\rho)\log\log(\rho)}{\log(R)-\log(\rho)}.$$
For every fixed $r_0>\rho$, one sees that
\begin{equation}\label{eta-1}
\eta_{1,R}(r_0)\sim\log\log(R)\quad\mbox{as }\:R\to\infty.
\end{equation}
For $k\ge 2$ the solutions to \eqref{ode-eta} can be represented as
$$\eta_{k,R}(r)=
A_{k,R}r^{\alpha_k^-}+B_{k,R}r^{\alpha_k^+}+\eta_k(r),$$
where
$$A_{k,R}=-\frac{R^{\alpha_k^+}\eta_k(\rho)-\rho^{\alpha_k^+}\eta_k(R)}
{R^{\alpha_k^+}\rho^{\alpha_k^-}-R^{\alpha_k^-}\rho^{\alpha_k^+}},\quad
B_{k,R}=-\frac{R^{\alpha_k^-}\eta_k(\rho)-\rho^{\alpha_k^-}\eta_k(R)}
{R^{\alpha_k^-}\rho^{\alpha_k^+}-R^{\alpha_k^+}\rho^{\alpha_k^-}},$$
and
\begin{equation}\label{e-eta}
\eta_k(r):=\frac{r^{\alpha_\ast}}{2\delta_k}
\left(
r^{\delta_k}\int_r^\infty\frac{t^{-\delta_k-1}}{\log^\sigma(t)}\,dt
+r^{-\delta_k}\int^r_\rho\frac{t^{\delta_k-1}}{\log^\sigma(t)}\,dt\right).
\end{equation}
It is easy to see that
\begin{equation}\label{eta-k}
0<\eta_k\le \frac{r^{\alpha_\ast}}{\delta_k^2\log^\sigma(r)},\qquad\forall k\ge 2.
\end{equation}
Moreover,
$\eta_k(r)=O(r^{\alpha_\ast}\log^{-\sigma}(r))$ as $r\to\infty$.

Represent $\psi=\sum_{k=1}^\infty \psi_k\phi_k$ as in \eqref{Fourier}, and set
\begin{equation}\label{e-log-sigma-infty}
w_{\psi,R}=\sum_{k=1}^\infty \eta_{k,R}\psi_k\phi_k.
\end{equation}
It is easy to see that the series converges in $H^1_0(\C_\Omega^{(\rho,R)})$
and $w_{\psi,R}$ solves \eqref{e-loglog-sum}.
Fix a compact $K\subset \C_\Omega^{\rho}$.
By Lemma \ref{l-series} and in view of \eqref{eta-k}, we conclude that
$$\sup_K\left|\sum_{k=2}^\infty \eta_{k,R}\psi_k\phi_k\right|\le
c\,\sup_K\frac{r^{\alpha_\ast}}{\delta_k^2\log^\sigma(r)}\le c_1,
\qquad\forall R\gg\rho,$$
with constants $c,c_1>0$ that do not depend on $R$.
Therefore
$$\inf_K w_{\psi,R}\sim \inf_K \left\{\eta_{1,R}\psi_1\phi_1\right\}\sim \log\log(R)\to\infty
\quad\mbox{as}\quad R\to\infty,$$
by \eqref{eta-1}, and the assertion follows.
\end{proof}

\begin{remark}
The value of $\sigma=1$ in the above Lemma is sharp.
When $\sigma>1$ it is not difficult to construct solutions to \eqref{e-log-sigma}
in the form \eqref{e-log-sigma-infty} with
$$\eta_1(r):=\frac{r^{\alpha_\ast}\log^{2-\sigma}(r)}{-(\sigma^2-3\sigma+2)}$$
and $\eta_k(r)$ as in \eqref{e-eta}.
Alternatively, if $\sigma>3/2$ then Hardy's inequality \eqref{e:Hardy} implies
that the quadratic form corresponding to \eqref{e-log-sigma}
satisfies the $\lambda$--property with $\lambda(x)=\frac{1/4}{|x|^2\log^2|x|}$.
Further, $|x|^{\alpha_\ast-2}\log^{-\sigma}|x|\in L^2(\lambda^{-1}dx)$ for any $\sigma>3/2$.
Thus Lemmas \ref{l:Riesz} and \ref{l:WMP} imply that \eqref{e-log-sigma-1}
has a positive solution in $D^1_0(\C_\Omega^\rho)$.
\end{remark}

\section{Proof of Theorem \ref{t:main'}, superlinear case $p\ge 1$}

We consider separately the cases $\mu<C_H+\lambda_1$ and $\mu=C_H+\lambda_1$.

\subsection{Case $\mu<C_H+\lambda_1$}

\paragraph{Nonexistence.}
First we prove the nonexistence of super-solutions in the subcritical case,
i.e. for $(p,s)$ below the critical line $\Lambda^\ast$.

\begin{lemma}\label{l:NONp>1<}
Let $p\ge 1$ and $s<\alpha^-(p-1)+2$.
Then \eqref{e:main'} has no positive super-solutions in $\C_\Omega^1$.
\end{lemma}
\begin{proof}
Let $w>0$ be a super-solution to (\ref{e:main'}) in $\C_\Omega^1$.
Then $w$ is a super-solution to the linear equation
\begin{equation}\label{e:01}
-\Delta w-\frac{\mu}{|x|^2}w=0\quad \mbox{in}\quad\C^1_\Omega.
\end{equation}
Choose a proper subdomain $\Omega^{\prime}\Subset\Omega$. Then, by
Theorem \ref{main:bound}, there exists $c>0$ such that
$$m_w(R,\Omega^{\prime})\ge cR^{\alpha^-}\quad(R\ge 2).$$
Linearizing \eqref{e:main'} and using the bound above, we conclude that $w>0$ is a super-solution to
\begin{equation}\label{e:06'}
-\Delta w-\frac{\mu}{|x|^2}w-\frac{V(x)}{|x|^2}w=0\quad\mbox{in
}\quad \C_{\Omega^{\prime}}^2,
\end{equation}
where $V(x):=Cw^{p-1}|x|^{2-s}$ satisfies
$$V(x)\geq c^{p-1}|x|^{\alpha^-(p-1)+(2-s)}\quad \mbox{in}\quad \C^2_{\Omega^{\prime}}.$$
Then the assertion follows by Corollary \ref{l:Nonexist}.
\end{proof}

Now we consider the critical case when $(p,s)$ belongs to the critical line $\Lambda^\ast$.

\begin{lemma}\label{l:NONp>1=}
Let $p\ge 1$ and $s=\alpha^-(p-1)+2$.
Then \eqref{e:main'} has no positive super-solutions in $\C_\Omega^1$.
\end{lemma}
\begin{proof}
Let $w>0$ be a super-solution to (\ref{e:main'}) in $\C_\Omega^1$.
Choose a proper subdomain $\Omega^{\prime}\Subset\Omega$. Arguing
as in the proof above we conclude that $w$ is a super-solution to
(\ref{e:06'}) with $V(x):=Cw^{p-1}|x|^{2-s}\geq\delta$ in
$\C^2_{\Omega^{\prime}}$, for some $\delta>0$.
Thus $w$ is a super-solution to the linear equation
\begin{equation}\label{e:07}
-\Delta w-\frac{W(\omega)}{|x|^2}w=0\quad\mbox{in}\quad\C^2_\Omega,
\end{equation}
where $W(\omega):=\mu+\eps\chi_{\Omega^\prime}$, with a fixed
$\eps\in(0,\delta]$. By the variational characterization of the
principal Dirichlet eigenvalue of
$-\Delta_\omega-\mu-\eps\chi_{\Omega^{\prime}}$ on $\Omega$ and
since $\mu<C_H+\lambda_1$, one can choose a small $\eps>0$
so that $C_H+\lambda_{W,1}>0$. Applying Theorem
\ref{main:bound} to (\ref{e:07}) we conclude that
$$m_w(R,\Omega^{\prime})\geq c r^{\alpha_{W,1}^-},\quad R\ge 4,$$
with $\alpha_1^-<\alpha_{W,1}^-<\alpha_\ast$.
Therefore $w>0$ is a super-solution to
$$-\Delta w-\frac{\mu}{|x|^2}w-\frac{\tilde V(x)}{|x|^2}w=0
\quad\mbox{in}\quad \C_{\Omega^{\prime}}^4,$$ where $\tilde
V(x):=Cw^{p-1}|x|^{2-s}$. Therefore
$$\tilde V(x)\geq C c^{p-1}|x|^{\alpha_{W,1}^- (p-1)+(2-s)}
\quad\mbox{in}\quad C^{4}_{\Omega^\prime},$$
with $\alpha_{W,1}^-(p-1)+(2-s)>0$.
Then the assertion follows from Corollary \ref{l:Nonexist}.
\end{proof}

\paragraph{Existence.}
Let $p>1$ and $s>\alpha^-_1(p-1)+2$.
Choose $\alpha\in(\alpha^-_1,\alpha^+_1)$ such that $\alpha\le\frac{s-2}{p-1}$.
Then one can verify directly that the functions
$$w:=\tau r^{\alpha}\phi_1(\omega)$$
are super-solution to \eqref{e:main'} in $\C_\Omega^1$ for sufficiently small $\tau>0$.
%Note that when $s<\alpha^+_1(p-1)+2$ and $\alpha=\frac{s-2}{p-1}$,
%then the super-solution $w$ could be extended to the entire cone $\C_\Omega$.
In the case $p=1$ and $s>2$ one sees that for any $\alpha\in(\alpha^-_1,\alpha^+_1)$
the function $w$ is a super-solution to \eqref{e:main'} in $\C_\Omega^\rho$ with a
sufficiently large $\rho\gg 1$.

\subsection{Case $\mu=C_H+\lambda_1$}

\paragraph{Nonexistence.}
The proof can be performed in one step for both subcritical and critical cases.

\begin{lemma}\label{l:NON_CH_p>1<}
Let $p\ge 1$ and $s\le\alpha_\ast(p-1)+2$.
Then \eqref{e:main'} has no positive super-solutions in $\C_\Omega^1$.
\end{lemma}
\begin{proof}
Assume that $w>0$ is a super-solution to (\ref{e:main'}) in $\C_\Omega^1$.
Then $w$ is a super-solution to
$$-\Delta w-\frac{C_H+\lambda_1}{|x|^2}w=0\quad \mbox{in}\quad\C^1_\Omega.$$
Choose a proper subdomain $\Omega^{\prime}\Subset\Omega$.
Then by Theorem \ref{main:bound}
$$m_w(R,\Omega^{\prime})\geq cR^{\alpha_\ast},\quad R\ge 2.$$
Linearizing \eqref{e:main'} and using the bound above, we conclude that $w>0$ is a super-solution to
\begin{equation}\label{e:06''}
-\Delta
w-\frac{C_H+\lambda_1}{|x|^2}w-\frac{W(x)}{|x|^2}w=0\quad\mbox{in
}\quad \C_{\Omega^{\prime}}^2,
\end{equation}
where $W(x):=Cw^{p-1}|x|^{2-s}\ge \tilde c$ in
$C_{\Omega^{\prime}}^2.$
%$$W(x)\geq c^{p-1}|x|^{\frac{(2-N)(p-1)}{2}+(2-s)}\quad \mbox{in}\quad C_{\Omega^{\prime}}^2,$$
%with $(2-N)(p-1)/2+(2-s)\geq 0$.
Then the assertion follows from Corollary \ref{l:Nonexist}.
\end{proof}

\paragraph{Existence.}
Let $p>1$ and $s>\alpha_\ast(p-1)+2$.
Choose $\beta\in(0,1)$.
Then one verifies directly that the functions
$$w:=\tau r^{\alpha_\ast}\log^\beta(r)\phi_1(\omega)$$
are super-solution to \eqref{e:main'} in $\C_\Omega^1$ for sufficiently small $\tau>0$.
In the case $p=1$ and $s>2$ one has to choose $\rho\gg 1$ sufficiently large.

\section{Proof of Theorem \ref{t:main'}, sublinear case $p<1$}

As before, we consider separately the cases $\mu<C_H+\lambda_1$ and $\mu=C_H+\lambda_1$.
First, we sketch the proofs of two auxiliary lemmas.

\begin{lemma}\label{l:EstP<1}
Let $p<1$. Let $w>0$ be a super-solution to \eqref{e:main'} in
$\C_\Omega^1$. Then for each proper subdomain
$\Omega^{\prime}\Subset\Omega$ there exists $c>0$ such that
\begin{equation}\label{sub}
m_w(R,\Omega^{\prime})\geq c\,R^{\frac{2-s}{1-p}},\qquad R\gg 1.
\end{equation}
\end{lemma}

\begin{proof}
Let $w>0$ be a super-solution to \eqref{e:main'}.
Then $-\Delta w\ge 0$ in $\C_\Omega^1$ and,
by the weak Harnack inequality (see, e.g. \cite[Theorem 8.18]{Gilbarg}),
for any $s>0$ and for any compact $K\subset\C_\Omega^1$ there exists $c>0$ such that
\begin{equation}\label{e:wHarnack}
\sup_{K} w^{-1}\leq c\left(R^{-N}\int_{K}w^{-s}dx\right)^{1/s}.
\end{equation}
%Clearly, $\frac{w^{p-1}}{|x|^s}\in L^1_{loc}(\C^1_\Omega)$.
Further, it follows from  Lemma \ref{l:A-ground} that
\begin{equation}\label{e:33}
\int_{\C^1_\Omega}|\nabla\varphi|^2dx-\mu\int_{\C^1_\Omega}\frac{\varphi^2}{|x|^2}dx\geq
C\int_{\C_\Omega^1}\frac{w^{p-1}}{|x|^s}\,\varphi^2dx,
\qquad\forall\,\varphi\in H^1_c(\C^1_\Omega)\cap H^\infty_c(\C^1_\Omega).
\end{equation}
Fix a proper subdomain $\Omega^\prime\Subset\Omega$.
Choose $\psi\in C_c^\infty(\Omega)$ such that $\psi=1$ on $\Omega^\prime$.
Choose $\theta_R(r)\in C^{0,1}_c(1,+\infty)$ such that
$0\leq\theta_R\leq 1$, $\theta_R=1$ for $r\in[R/2,R]$,
$Supp(\theta_R)=[R/4,2R]$ and $|\nabla\theta_R|<c/R$.
Then
\begin{equation}\label{e:31}
\int_{\C^1_\Omega}|\nabla(\theta_R\psi)|^2dx-\mu\int_{\C^1_\Omega}\frac{|\theta_R\psi|^2}{|x|^2}dx\leq cR^{N-2}.
\end{equation}
%\begin{eqnarray}\label{e:31}
%&&\int_{\C^1_\Omega}|\nabla\eta|^2dx-\mu\int_{\C^1_\Omega}\frac{\eta^2}{|x|^2}dx\leq
%\int_{\C^1_\Omega}\left(|\nabla_r\theta|^2\phi_1^2+\theta^2|\nabla_\omega\phi_1|^2\right)\\&&=
%\int_{R/2}^{R}|\nabla_r\theta|^2
%r^{N-1}dr+\lambda_1(\Omega')\int_{R/4}^{2R} \theta^2r^{N-1} dr\leq
%cR^{N-2}.\nonumber
%\end{eqnarray}
On the other hand,
\begin{equation}\label{e:32}
\int_{\C^1_{\Omega^\prime}}\frac{w^{p-1}}{|x|^s}\,(\theta_R\psi)^2\,dx\geq
\int_{\C^{(R/2,R)}_{\Omega^\prime}}\frac{w^{p-1}}{|x|^s}dx\geq
R^{-s}\int_{\C^{(R/2,R)}_{\Omega^\prime}}w^{p-1}dx.
\end{equation}
Combining \eqref{e:33}, \eqref{e:31} and \eqref{e:32} we derive
$$cR^{s-2}\geq R^{-N}\int_{\C^{(R/2,R)}_{\Omega^\prime}}w^{p-1}dx.$$
%Applying Harnack inequality
%
By \eqref{e:wHarnack} we obtain
$$cR^{\frac{s-2}{1-p}}\geq
\left(R^{-N}\int_{\C^{(R/2,R)}_{\Omega^\prime}}w^{-(1-p)}\,dx\right)^{\frac{1}{1-p}}\geq
c_1\sup_{\C_{\Omega^\prime}^{(R/2,R)}}w^{-1}.$$ Hence the assertion
follows.
\end{proof}

\begin{lemma}\label{l:sub-super}
Let $p<1$, $\mu\le C_H+\lambda_1$ and $s\in\R$. Assume that
\eqref{e:main'} has a positive super-solution in $\C_\Omega^1$.
Then
%\eqref{e:main'} has a positive solution in $\C_\Omega^1$.
there exists a positive solution to \eqref{e:main'} in $\C_\Omega^1$.
\end{lemma}

\begin{proof}
Let $w>0$ be a super-solution to \eqref{e:main'} in $\C_\Omega^1$.
Then, by Lemma \ref{l:Minimal} or \ref{l:Minimal-eps}, $w\ge
cv_\psi$ in $\C_\Omega^1$, where $v_\psi$ is a comparison function
defined by \eqref{e:v-psi} or \eqref{e:v-psi-eps}. Obviously,
$v_\psi>0$ is a sub-solution to \eqref{e:main'} in $\C_\Omega^2$.
Thus we can proceed via the standard sub and super-solutions techniques
to prove existence of a solution to \eqref{e:main'} in $\C_\Omega^2$,
located between $cv_\psi$ and $w$ (cf. \cite[Proposition 1.1(iii)]{KLM}).
Finally, after a suitable scaling we obtain a solution to \eqref{e:main'} in $\C_\Omega^1$.
\end{proof}

\subsection{Case $\mu<C_H+\lambda_1$}

\paragraph{Nonexistence.}
We distinguish between the subcritical and critical cases.
When $(p,s)$ is below the critical
line, the proof of the nonexistence is straightforward.

\begin{lemma}\label{l:NONp<1<}
Let $p<1$ and $s<\alpha_1^+(p-1)+2$.
Then \eqref{e:main'} has no positive super-solutions in $\C_\Omega^1$.
\end{lemma}
\begin{proof}
Let $w>0$ be a super-solution to (\ref{e:main'}).
Then $w$ is a super-solution to the linear equation
\begin{equation}\label{e:01'}
-\Delta w - \frac{\mu}{|x|^2}w=0\quad\mbox{in }\:\C_\Omega^1\, .
\end{equation}
Choose a proper subdomain $\Omega^{\prime}\Subset\Omega$. By
Theorem \ref{main:bound} we conclude that
\begin{equation}\label{e:u}
m_w(R,\Omega^{\prime})\le c R^{\alpha_1^+},\qquad R\gg 1.
\end{equation}
This contradicts to \eqref{sub}.
\end{proof}

Next we consider the case when $(p,s)$ is on the critical line,
and hence \eqref{e:u} is no longer incompatible with \eqref{sub}.

\begin{lemma}\label{l:NONp<1=}
Let $p<1$ and $s=\alpha_1^+(p-1)+2$.
Then \eqref{e:main'} has no positive super-solutions in $\C_\Omega^1$.
\end{lemma}

\begin{proof}
Let $w>0$ be a super-solution to \eqref{e:main'}. According to
Lemma \ref{l:sub-super} we may assume that $w$ is a solution to
\eqref{e:main'}. Choose a proper subdomain
$\Omega^{\prime}\Subset\Omega$. Linearizing \eqref{e:main'} and
using the upper bound \eqref{sub} we conclude that $w>0$ is a solution to
\begin{equation}\label{e:06'''}
-\Delta w-\frac{\mu}{|x|^2}w-\frac{V(x)}{|x|^2}w=0\quad\mbox{in }\:\C_{\Omega}^1,
\end{equation}
where $V(x):=c^{p-1}|x|^{2-s}w^{p-1}$ satisfies $V(x)\leq c_1$ in
$\C_{\Omega^{\prime}}^{\rho_1}$, with a fixed $\rho_1\gg 1$.
This implies, in particular, that $w$ satisfies strong Harnack's inequality with
$r$--independent constants. More precisely, for a given subdomain
$\Omega^{\prime\prime}\Subset\Omega^{\prime}$ one has
\begin{equation}\label{e:u-strong-harnack}
M_w(R,\Omega^{\prime\prime})\le
C_s\,m_w(R,\Omega^{\prime\prime}),\qquad R\gg\rho,
\end{equation}
where $C_s=C_s(\Omega^{\prime\prime})>0$ does not depend on $R\gg\rho$.
Using \eqref{e:u-strong-harnack} and the upper bound \eqref{e:u} we conclude that
\begin{equation}\label{e:u-strong}
M_w(R,\Omega^{\prime\prime})\le c_2R^{\alpha_1^+},\qquad R\gg\rho.
\end{equation}
This implies that $V(x)\ge\delta$ in $\C_{\Omega^{\prime\prime}}^{\rho_2}$,
for some $\delta>0$ and $\rho_2\gg\rho_1$.
Hence $w>0$ is a super-solution to the linear equation
\begin{equation}\label{lin-3}
-\Delta w-\frac{W_\varepsilon(\omega)}{|x|^2}w=0\quad\mbox{in
}\C_{\Omega}^{\rho_2},
\end{equation}
where $W_\eps(\omega):=\mu+\eps\chi_{\Omega^{\prime\prime}}$,
with a fixed $\eps\in(0,\delta]$.
By the variational characterization of the principal Dirichlet eigenvalue of
$-\Delta_\omega-W_\eps$ on $\Omega$ and since
$\mu<C_H+\lambda_1$, one can choose a small $\eps>0$ such that
$C_H+\lambda_{W_\eps,1}>0$.
Applying Theorem \ref{main:bound} to \eqref{lin-3} we conclude that
\begin{equation}\label{e:u-improved}
m_w(R,\Omega^{\prime})\le c_2R^{\alpha_{W_\eps,1}^+}.
\end{equation}
with $\alpha_{W_\eps,1}^+<\alpha_1^+$.
Now \eqref{e:u-improved} contradicts to upper bound \eqref{sub}.
\end{proof}

\paragraph{Existence.}
Let $s>\alpha^+_1(p-1)+2$.
Assume that $0\le p<1$.
Choose $\alpha\in(\alpha^-_1,\alpha^+_1)$ such that $\alpha\ge\frac{s-2}{p-1}$.
Then there exists a unique bounded positive solution to the problem
$$-\Delta_\omega\phi-(\alpha(\alpha+N-2)+\mu)\phi=1,\qquad \phi\in H^1_0(\Omega).$$
Further, a direct computation verifies that the functions
$$w:=\tau r^{\alpha}\phi(\omega)$$
are super-solutions to \eqref{e:main'} in $\C_\Omega^1$ for a sufficiently large $\tau>0$.

Now assume that $p<0$. Choose $\alpha$ as above,
so there exists a unique bounded positive solution of the problem
$$-\Delta_\omega\bar\phi-(\alpha(\alpha+N-2)+\mu)(\bar\phi+1)=1,\qquad \bar\phi\in H^1_0(\Omega).$$
%Choose $\tau>0$ such that $c\tau^{p-1}\|\bar\phi\|_{L^\infty}^p=1.$
Then
$$w:=\tau r^{\alpha}\bar\phi(\omega)$$
are super-solutions to \eqref{e:main'} in $\C_\Omega^1$ for sufficiently large $\tau>0$.

\subsection{Case $\mu=C_H+\lambda_1$}

\paragraph{Nonexistence.}
The proof is straightforward for $(p,s)$  below the critical
line $\Lambda$.

\begin{lemma}\label{l:NON CH p<1<}
Let $p<1$ and $s<\alpha_\ast(p-1)+2$.
Then \eqref{e:main'} has no positive super-solutions in $\C_\Omega^1$.
\end{lemma}
\begin{proof}
Let $w>0$ be a super-solution to (\ref{e:main'}).
Similarly to the proof of Lemma \ref{l:NONp<1<},
by Theorem \ref{main:bound} we conclude that for a proper subdomain $\Omega^\prime\Subset\Omega$,
\begin{equation}\label{e:u-log}
m_w(R,\Omega^{\prime})\le c R^{\frac{2-N}{2}}\log(R), \qquad R\gg 1.
\end{equation}
This contradicts to the upper bound \eqref{sub}.
\end{proof}

When $(p,s)$ belongs to the critical line $\Lambda$ inequality
\eqref{e:u-log} is no longer incompatible with \eqref{sub}.

\begin{lemma}\label{l:NON_CH_p=-1=}
Let $p\in[-1,1)$ and $s=\alpha_\ast(p-1)+2$.
%Let $p=-1$. Assume that $\mu=C_H+\lambda_1$ and $s=N$.
Then \eqref{e:main'} has no positive super-solutions in $\C_\Omega^1$.
\end{lemma}
\begin{proof}
Let $w>0$ be a super-solution to \eqref{e:main'}.
If $p\in[0,1)$ then the lower bound \eqref{sub} implies that $w$ is a super-solution to
\begin{equation}\label{lin-5}
-\Delta w-\frac{C_H+\lambda_1}{|x|^2}w=\frac{\psi(\omega)}{|x|^{2-\alpha_\ast}}
\quad\mbox{in}\quad\C_{\Omega}^{\rho^\prime},
\end{equation}
with some $\psi\in C^\infty_c(\Omega)$ and $\rho^\prime>\rho$.
From Lemma \ref{l-log-sigma}, it follows that \eqref{lin-5} has no positive super-solutions.

Let $p\in[-1,0)$.
According to Lemma \ref{l:sub-super} we may assume that $w$ is a solution to \eqref{e:main'}.
Similarly to the proof of Lemma \ref{l:NONp<1=}, for a proper subdomain $\Omega^{\prime}\Subset\Omega$ and
a function $\psi\in C^\infty_c(\Omega^\prime)$
we conclude that $w>0$ is a super-solution to the linear equation
\begin{equation}\label{lin-6}
-\Delta w-\frac{C_H+\lambda_1}{|x|^2}w=\frac{\psi(\omega)}{|x|^{2-\alpha_\ast}\log^{-p}|x|}
\quad\mbox{in}\quad\C_{\Omega}^{\rho^\prime},
\end{equation}
for some $\rho^\prime>\rho$.
Then the assertion follows from Lemma \ref{l-log-sigma}.
\end{proof}

\paragraph{Existence.}
In the critical case $\mu=C_H+\lambda_1$ positive super-solutions to \eqref{e:main'} with $p<1$
can not be constructed as "pseudo"--radial functions of the form $u=v(r)\varphi(\omega)>0$,
as the following proposition shows.

\begin{proposition} %Pseudo--radial super-solutions to $\omega$--critical operator
Let $u=v(r)\varphi(\omega)>0$ be a super-solution to
\begin{equation}\label{e-CH}
-\Delta u-\frac{C_H+\lambda_1}{|x|^2}u=0\quad\mbox{in}\quad\C_\Omega^\rho.
\end{equation}
Then $u=v(r)\phi_1(\omega)$, where $v$ is a super-solution to
\begin{equation}\label{e-CH-rad}
-\frac{\partial^2 v}{\partial r^2}-\frac{N-1}{r}\frac{\partial v}{\partial r}
-\frac{C_H}{r^2}v\ge 0\quad\mbox{in}\quad(\rho,\infty).
\end{equation}
\end{proposition}
\begin{proof}
Let $u=v(r)\varphi(\omega)>0$ be a super-solution to \eqref{e-CH}. Then
$$\left\{-\frac{\partial^2 v}{\partial r^2}-\frac{N-1}{r}\frac{\partial v}{\partial r}
-\frac{C_H}{r^2}v\right\}\varphi+\left\{(-\Delta_\omega-\lambda_1)\varphi\right\}\frac{v}{r^2}\ge 0
\quad\mbox{in}\quad\C_\Omega^\rho.$$
Separating the variables and using Barta's inequality (see Lemma \ref{l:A-ground})
$$\sup_{0<\varphi\in H^1_{loc}(\Omega)}
\left\{\inf_{\omega\in\Omega}\frac{(-\Delta_\omega-\lambda_1)\varphi}{\varphi}\right\}\le 0,$$
we obtain
$$\frac{r^2}{v}\left\{-\frac{\partial^2 v}{\partial r^2}-\frac{N-1}{r}\frac{\partial v}{\partial r}
-\frac{C_H}{r^2}v\right\}\ge
-\left\{\inf_{\omega\in\Omega}\frac{(-\Delta_\omega-\lambda_1)\varphi}{\varphi}\right\}\ge
0\quad\mbox{in}\quad\C_\Omega^\rho.$$
On the other hand, the one--dimensional Hardy's inequality implies that the inequality
$$-\frac{\partial^2 v}{\partial r^2}-\frac{N-1}{r}\frac{\partial v}{\partial r}
-\frac{C_H}{r^2}v\ge\epsilon\frac{v}{r^2}\quad\mbox{in}\quad(\rho,\infty)$$
has a positive solution if and only if $\epsilon=0$. Hence
$$\inf_{\omega\in\Omega}\frac{(-\Delta_\omega-\lambda_1)\varphi}{\varphi}=0,$$
and therefore, $\varphi=\phi_1$. We conclude
that $u$ must be of the form $u=v(r)\phi_1(\omega)$, where $v$ is a super-solution to \eqref{e-CH-rad}.
\end{proof}

It is easy to see that if $\Omega\Subset S^{N-1}$ is a proper subdomain of the sphere
then equation \eqref{e:main'} with $p<1$ does not admit positive super-solutions of the form
$v(r)\phi_1(\omega)$. Nevertheless, for $(p,s)$ above the critical line $\Lambda$ we prove the following.

\begin{lemma}
Let $p<1$ and $s>\alpha_\ast(p-1)+2$.
Then
%\eqref{e:main'} has a positive super-solution in $\C_\Omega^1$.
there exists a positive super-solution to \eqref{e:main'}.
\end{lemma}
\begin{proof}
Given $\eps\in(0,1/4)$, $\sigma>3/2$ and $\rho\ge\exp(1)$, consider the problem
\begin{equation}\label{e-log-sigma-1}
-\Delta w-\frac{C_H+\lambda_1}{|x|^2}w-\frac{\eps}{|x|^2\log^2|x|}w=\frac{1}{|x|^{2-\alpha_\ast}\log^\sigma|x|},
\qquad w\in\D(\E_V,\C_\Omega^\rho).
\end{equation}
It follows from Hardy's inequality \eqref{e:Hardy} that the quadratic form that corresponds to \eqref{e-log-sigma-1}
satisfies the $\lambda$--property with $\lambda(x)=\frac{1/4-\eps}{|x|^2\log^2|x|}$.
Further, $|x|^{\alpha_\ast-2}\log^{-\sigma}|x|\in L^2(\lambda^{-1}dx)$.
Thus Lemmas \ref{l:Riesz} and \ref{l:WMP} imply that \eqref{e-log-sigma-1} has a unique solution $w_\sigma>0$.
Choose $\beta>\sigma$ and set
$$v_{\sigma}:=w_\sigma+\frac{|x|^{\alpha_\ast}}{\log^\beta|x|}.$$
Then
$$\left(-\Delta-\frac{C_H+\lambda_1}{|x|^2}-\frac{\eps}{|x|^2\log^2|x|}\right)v_{\sigma}=
\frac{1}{|x|^{2-\alpha_\ast}}
\left(\frac{1}{\log^\sigma|x|}-\frac{\lambda_1}{\log^{\beta}|x|}-\frac{\beta(1+\beta)+\eps}{\log^{\beta+2}|x|}\right)
\ge 0\mbox{ in }\C_\Omega^{\rho^\prime},$$
for some $\rho^\prime\gg\rho$.
Set $\delta:=s-\alpha_\ast(p-1)-2>0$ and choose $\tau=\tau(\delta)>0$ such that
$$\frac{C}{|x|^s}(\tau v_\sigma)^{p-1}\leq
\frac{C\tau^{p-1}\log^{\beta(1-p)+2}|x|}{\eps|x|^{\delta}}\frac{\eps}{|x|^2\log^2|x|}\le\frac{\eps}{|x|^2\log^2|x|}
\quad\mbox{in}\quad\C_\Omega^{\rho^\prime}.$$
Then
$$\left(-\Delta-\frac{C_H+\lambda_1}{|x|^2}\right)(\tau v_\sigma)\ge
\frac{\eps}{|x|^2\log^2|x|}(\tau v_\sigma)\ge
\frac{C}{|x|^s}(\tau v_\sigma)^{p-1}(\tau v_\sigma)=\frac{C}{|x|^s}(\tau v_\sigma)^{p}
\quad\mbox{in}\quad\C_\Omega^{\rho^\prime},$$
that is $\tau v_\sigma>0$ is a super-solution to \eqref{e:main'} in $\C_\Omega^{\rho^\prime}$.
\end{proof}

In the case of exterior domains,
the existence of positive super-solutions on the critical line $\Lambda$ for $p<-1$ is elementary observed.

\begin{lemma}
Assume that $\Omega=S^{N-1}$.
Let $p<-1$ and $s=\alpha_\ast(p-1)+2$.
Choose $\beta\in(\frac{2}{1-p},1)$.
Then
$$v:=\tau|x|^{\alpha_\ast}\log^\beta|x|$$
is a super-solution to \eqref{e:main'} in $\C_\Omega^{\rho}$ for sufficiently large $\tau>0$.
\end{lemma}

\begin{proof}
The direct computation.
\end{proof}

In the case of proper domains $\Omega\Subset S^{N-1}$,
the existence (or nonexistence) of positive super-solutions to \eqref{e:main'}
with $p<-1$ and $s=\alpha_\ast(p-1)+2$
becomes a more delicate issue that remains open at the moment.
The analysis of the decay rate of super-solutions to \eqref{e-CH}
near the lateral boundary of the cone should be invoked.
%Notice that after the transformation $u(r,\omega)=r^{\alpha_\ast}\,v$,
%where $v=v(t,\omega)$, $t=t(r)=-\alpha_\ast\log(r)$,
%and a scaling, equation \eqref{e:main'} with $\mu=C_H+\lambda_1$ and $s=\alpha_\ast(p-1)+2$
%%transforms into a particulary simple form equation
%takes a particularly simple form
%$$v_{tt}+(-\Delta_\omega-\lambda_1) v=C\,v^p\quad
%\mbox{in}\quad(\log(\rho),\infty)\times\Omega.$$
We will return to this problem elsewhere.

\appendix

\section{Appendix}\label{A1}

Let $\E_V$ be a symmetric bilinear form defined by
$$\E_V(u,v):=\int_G\nabla u\cdot\nabla v\,dx-\int_G Vuv\,dx\qquad(u,v\in H^1_c(G)\cap L^\infty_c(G)),$$
where $G\subseteq\R^N$ is a domain and $0\leq V\in L^1_{loc}(G)$ a potential.
Below we present several facts concerning the relations between the positivity
of the form $\E_V$ and the existence of positive (super)\,solutions
to the linear equation
\begin{equation}\label{e:linear-V}
(-\Delta-V)v=f\quad \mbox{in}\quad G,
\end{equation}
associated with $\E_V$, where $f\in L^1_{loc}(G)$.
A super-slolution to \eqref{e:linear-V} is a function $u\in H^1_{loc}(G)\cap L^1_{loc}(G,Vdx)$ such that
\begin{equation}\label{e:linear-weak}
\int_G\nabla u\cdot\nabla\varphi\,dx-\int_G Vu\varphi\,dx\ge\int_G f\varphi\,dx,
\qquad\forall\:0\le\varphi\in H^1_c(G)\cap L^\infty_c(G).
\end{equation}
The notions of a sub-solution and solution are defined similarly
by replacing "$\ge$" with "$\le$" and "$=$" respectively.
Most of the results below
are known from the theory of Dirichlet forms (cf. \cite{Davies,FOT})
and Agmon's criticality theory (cf. \cite{Agmon,Agmon-2}).
We include the proofs for the completeness and reader's convenience.

\paragraph{Extended Dirichlet Space.}
Assume that the form $\E_V$ is {\sl positive definite}, that is
\begin{equation}\label{e:PD}
\E_V(u,u)>0,\qquad\:\forall\:0\neq u\in H^1_c(G)\cap L^\infty_c(G).
\end{equation}
Following Fukushima \cite[p.35--36]{FOT},
denote by $\D(\E_V,G)$ the family of measurable a.e.\ finite functions $u:G\to\R$
such that there exists an $\E_V$--Cauchy sequence $(u_n)\subset H^1_c(G)\cap L^\infty_c(G)$
that converges to $u$ a.e. in $G$. This sequence $(u_n)$ is called an {\sl approximating sequence}
for $u\in\D(\E_V,G)$. Then the limit $\E_V(u,u):=\lim_{n\to\infty}\E_V(u_n,u_n)$ exists
and is independent of the choice of the approximating sequence.
Thus $\E_V$ is extended uniquely to a nonnegative definite bilinear form on $\D(\E_V,G)$.
The family $\D(\E_V,G)$ is called the {\sl extended Dirichlet space} of $\E_V$.
It is not a Hilbert space, in general.
However, $\D(\E_V,G)$ is invariant under the standard truncations.

\begin{lemma}\label{l:pm}
Let $u\in\D(\E_V,G)$. Then $u^+=u\vee 0\in\D(\E_V,G)$, $u^-=-(u\wedge 0)\in\D(\E_V,G)$ and
\begin{equation}\label{e:pm}
\E_V(u^\pm,u^\pm)\le\E(u,u),\qquad\forall u\in\D(\E_V,G).
\end{equation}
If $u,v\in\D(\E_V,G)$ then $u\vee v,u\wedge v\in\D(\E_V,G)$.
\end{lemma}
\begin{proof}
Assume $u\in H^1_c(G)\cap L^\infty_c(G)$. Then $u^+\in
H^1_c(G)\cap L^\infty_c(G)$. By the direct computation we have
$$\E_V(u^+,u^+)+\E_V(u^-,u^-)=\E_V(u,u).$$
Hence \eqref{e:pm} follows by \eqref{e:PD} for any $u\in
H^1_c(G)\cap L^\infty_c(G)$, and, then, for arbitrary
$u\in\D(\E_V,G)$ by a standard approximation argument.
%Using \eqref{e:PD} we obtain
%\begin{eqnarray*}
%\E_V(u^+,u^+)&=&\int_G|\nabla u^+|^2\,dx-\int_G Vu^2\,dx+\int_G V(u^-)^2\,dx\\
%&\le&\int_G|\nabla u^+|\,dx-\int_G Vu^2\,dx+\int_G|\nabla u^-|^2\,dx=\E_V(u,u).
%\end{eqnarray*}
%Now \eqref{e:pm} follows for arbitrary $u\in\D(\E_V,G)$ by a standard approximation argument.
%Similarly for $u^-$. Hence the assertion for $u\vee v$ and $u\wedge v$ follows.
\end{proof}
\begin{remark}
We do not claim that $u\in\D(\E_V,G)$ implies $u\wedge 1\in\D(\E_V,G)$.
\end{remark}

Following \cite{Agmon-2,Allegretto}, we say that the form $\E_V$ satisfies
the {\sl $\lambda$--property} if there exists a function $0<\lambda\in L^1_{loc}(G)$ such that
$\lambda^{-1}\in L^\infty_{loc}(G)$ and
\begin{equation}\label{e:lambda}
\E_V(u,u)\ge\int_G u^2\,\lambda(x)\,dx,\qquad\forall\:u\in H^1_c(G)\cap L^\infty_c(G).
\end{equation}
If $\E_V$ satisfies the $\lambda$--property then the extended Dirichlet space $\D(\E_V,G)$
is a Hilbert space with the inner product $\E_V(\cdot,\cdot)$ and
the corresponding norm $\|\cdot\|_\D=\sqrt{\E_V(\cdot,\cdot)}$.
Clearly
$$H^1_{c}(G)\cap L^\infty_c(G)\subset\D(\E_V,G)\subset H^1_{loc}(G)
\quad\mbox{and}\quad
\D(\E_V,G)\subset L^2(G,\lambda\,dx).$$
By $\D^\prime(\E_V,G)$ we denote the space of linear continuous functionals on $\D(\E_V,G)$.
The following lemma is a standard consequence of the Riesz Representation Theorem.
\begin{lemma}\label{l:Riesz}
Assume that $\E_V$ satisfies the $\lambda$--property.
Let $l\in\D^\prime(\E_V,G)$.
Then there exists a unique $\phi_\ast\in\D(\E_V,G)$ such that
\begin{equation}\label{Riesz}
\E_V(\phi_\ast,\varphi)=l(\varphi),\qquad\forall\:\varphi\in\D(\E_V,G).
\end{equation}
\end{lemma}

Let $\hat\D^\prime(\E_V,G):=\{f\in L^1_{loc}(G):\int_G f\varphi\,dx\le c\|\varphi\|_\D\:,\:
\forall\varphi\in H^1_c(G)\cap L^\infty_c(G)\}$.
It is easy to see that
$$L^2(G,\lambda^{-1}\,dx)\subset\hat\D^\prime(\E_V,G).$$
Clearly $\hat\D^\prime(\E_V,G)$ can be identified with a linear subspace of $\D^\prime(\E_V,G)$.
Thus Lemma \ref{l:Riesz} implies that for any $f\in\hat\D^\prime(\E_V,G)$ the problem
\begin{equation}\label{e:Lax}
(-\Delta-V)u=f,\qquad u\in\D(\E_V,G),
\end{equation}
has a unique solution.

\paragraph{Maximum and comparison principles.}
Consider the homogeneous equation
\begin{equation}\label{e:ground}
(-\Delta-V)u=0\quad\mbox{in }\:G.
\end{equation}
We present weak maximum and comparison principles
for solutions and super-solutions of \eqref{e:linear-V} in a form suitable for our framework.

\begin{lemma}\label{l:WMP}
%{\sf (Weak Maximum Principle)}
Assume that $\E_V$ satisfies the $\lambda$--property.
Let $w\in H^1_{loc}(G)$ be a super-solution to \eqref{e:ground} such
that $w^-\in \D(\E_V,G)$. Then $w\geq 0$ in $G$.
\end{lemma}
\begin{proof}
Let $(\varphi_n)\subset H^1_{c}(G)\cap L^\infty_c(G)$ be an approximating sequence
for $w^-\in \D(\E_V,G)$. Set $w_n:=\varphi_n^+\wedge w^-$.
Hence $0\le w_n\in \D(\E_V,G)$, by Lemma \ref{l:pm}.
Note that $w_n=w^-+(\varphi_n^+-w^-)^-$.
Therefore
$$\E_V(w^--w_n,w^--w_n)=\E_V((\varphi_n^+-w^-)^-,(\varphi_n^+-w^-)^-)
\le\E_V(\varphi_n-w^-,\varphi_n-w^-)\to 0.$$
Thus $(w_n)$ is a nonnegative approximating sequence for $w^-$.
Since $w^+\wedge w_n=0$, we obtain
$$0\le\E_V(w,w_n)=-\E_V(w^-,w_n)\to-\E_V(w^-,w^-)\le 0.$$
We conclude that $w^-=0$.
\end{proof}

\begin{remark}
Note that if $u\gneq 0$ is a super-solution to \eqref{e:ground}
then $-\Delta u\ge 0$ in $G$.
Hence, by the by the weak Harnack inequality, $u>0$ in $G$.
\end{remark}

\begin{remark}\label{r:WMP}
If $\E_V$ satisfies the $\lambda$--property then Lemmas \ref{l:Riesz} and \ref{l:WMP}
imply that equation \eqref{e:ground} has a rich cone of positive super-solutions.
Indeed, if $0\lneq f\in \hat\D^\prime(\E_V,G)$ and
$\phi_\ast\in\D(\E_V,G)$ is the solution to \eqref{e:Lax}
then $\phi_\ast>0$ in $G$.
\end{remark}

The following comparison principle is a straightforward consequence of Lemma \ref{l:WMP}.

\begin{corollary}\label{c:A7}
Assume that $\E_V$ satisfies the $\lambda$--property.
Let $w\in H^1_{loc}(G)$ be a super-solution to \eqref{e:ground} and
$v\in H^1_{loc}(G)$ be a sub-solution to \eqref{e:ground} such that $(w-v)^-\in \D(\E_V,G)$.
Then $w\geq v$ in $G$.
\end{corollary}

A version of the comparison principle below plays a crucial role in
the analysis of asymptotic behavior of super-solutions to \eqref{e:V}
in Section \ref{s:asymptotic}.

\begin{lemma}\label{l:WCP}
%{\sf (Weak Comparison Principle)}
Assume that $\E_V$ satisfies the $\lambda$--property.
Let $0\leq w\in H^1_{loc}(G)$, $v\in \D(\E_V,G)$ and
$w-v$ be a super-solution to equation \eqref{e:ground}.
Then $w\geq v$ in $G$.
\end{lemma}
\proof
Let $(G_n)$ be an exhaustion of $G$,
i.e.\ an increasing sequence of bounded smooth domains such that
$G_n\Subset G_{n+1}\Subset G$ and $\cup G_n=G$.
Note that $\lambda^{-1}\in L^\infty(G_n)$ and therefore $\D(\E_V,G_n)=H^1_0(G_n)$.
Clearly $H^1_0(G_n)$ is a closed subspace of $\D(\E_V,G)$.

Let $v\in \D(\E_V,G)$.
Let $f\in \D^\prime(\E_V,G)$ be defined by
$$f(\varphi):=\E_V(v,\varphi),\qquad(\varphi\in \D(\E_V,G)).$$
By Lemma \ref{l:Riesz} there exists the unique $v_n\in H^1_0(G_n)$
such that
$$\E_V(v_n,\varphi) = f(\varphi),\qquad\forall\varphi\in H^1_0(G_n).$$
Thus
$$(-\Delta-V)(v-v_n)= 0\quad\mbox{in }\:G_n,$$
and hence
$$(-\Delta-V)(w-v_n)\ge 0\quad\mbox{in }\:G_n,$$
with $w-v_n\in H^1_{loc}(G_n)$ and $0\le(w-v_n)^-\le v_n^+\in
H^1_0(G_n)$. Corollary \ref{c:A7} implies  %We conclude
that $v_n\le w$ in $G_n$.

Let $\bar{v}_n$ denote the extension of $v_n$ to $G$ by zero.
Clearly $\bar{v}_n\in\D(\E_V,G)$.
To complete the proof it suffices to show that
$\bar{v}_n\to v$ in $\D(\E_V,G)$.
Indeed, by the construction of $\bar{v}_n$ we obtain
$$\E_V(\bar{v}_n,\bar{v}_n)=f(v_n)\le
\|f\|_{\D^\prime}\|\bar{v}_n\|_\D.$$
Hence the sequence $(\bar{v}_n)$ is bounded in $\D(\E_V,G)$.
Thus there is a subsequence, which we still denote by $(\bar{v}_n)$,
that converges weakly to $v_\ast\in \D(\E_V,G)$.
Now let $\varphi\in H^1_c(G)\cap L^\infty_c(G)$.
Then $\varphi\in H^1_0(G_n)$ for all $n\in\N$ large enough, and
$$\E_V(\bar{v}_n,\varphi)=f(\varphi).$$
Passing to the limit we conclude that
$$\E_V(v_\ast,\varphi)=f(\varphi),\qquad\forall\varphi\in H^1_c(G)\cap L^\infty_c(G),$$
and therefore $v_\ast=v$. Furthermore,
$$\E_V(\bar{v}_n-v,\bar{v}_n-v)=f(\bar{v}_n)-2f(v)+f(v).$$
Since $f(\bar{v}_n)\to f(v)$, it follows that $\bar{v}_n\to v$ strongly in $\D(\E_V,G)$.
\qed

\paragraph{Ground state transformation.}
If $\E_V$ satisfies the $\lambda$--property,
then equation \eqref{e:ground} has a rich cone of positive super-solutions, see Remark \ref{r:WMP}.
One can show that \eqref{e:ground} has a positive solution if $\E_V$ is positive definite
(but may not satisfy the $\lambda$--property, cf. \cite[Theorem 3.1]{Agmon}).
Below we prove the converse
(cf. \cite{Agmon}, \cite{Davies} for the ground state transform,
\cite{Allegretto} for the Picone identity,
\cite{Pinsky} for the $h$--transform).
\begin{lemma}\label{l:A-ground}
Let $0<\phi\in H_{loc}^1(G)$ be a (super)\,solution to the equation
\begin{equation}\label{e:ground-f}
(-\Delta-V)\phi=f\quad \mbox{in}\quad G,
\end{equation}
where $0\le f\in L^1_{loc}(G)$.
Then the form $\E_V$ is positive definite in the sense of \eqref{e:PD}.
Moreover,
\begin{equation}\label{ground}
\D(\E_V,G)\ni u\mapsto\frac{u}{\phi}\in H^1(G,\phi^2 dx)
\end{equation}
and
\begin{equation}\label{A-ground}
\E_V(u,u)\,(\ge)=\,\int_G|\nabla\left(\frac{u}{\phi}\right)|^2\phi^2 dx+\int_G u^2\,\frac{f}{\phi}\,dx,
\qquad\forall u\in \D(\E_V,G).
\end{equation}
\end{lemma}
\begin{proof}
Let $u\in H^1_c(G)\cap L^\infty_c(G)$.
Then $\varphi:=\frac{u^2}{\phi}\in H_{c}^{1}(G)\cap L^\infty_c(G)$.
Testing \eqref{e:ground} against $\varphi$ we arrive at
$$2\int_G u\,\nabla u\,\frac{\nabla \phi}{\phi}\,dx\,(\ge)=\,\int_G u^2 \frac{|\nabla
\phi|^2}{\phi^2}\,dx+\int_G Vu^2\,dx+\int_G \frac{f}{\phi}u^2\,dx.$$
Direct computation gives that
\begin{eqnarray*}&&
\int_G \left(|\nabla u|^2-Vu^2\right)dx-\int_G|\nabla\left(\frac{u}{\phi}\right)|^2\phi^2dx\\
&&\qquad\qquad=\int_G \left(|\nabla u|^2-Vu^2\right)dx-\int_G\left(
\frac{|\nabla u|^2}{\phi^2}-2u\,\nabla u\,\frac{\nabla \phi}{\phi^3}+
u^2\frac{|\nabla\phi|^2}{\phi^4} \right)\phi^2dx\\
&&\qquad\qquad=2\int_G u\,\nabla u\,\frac{\nabla \phi}{\phi}\,dx-\int_G u^2
\frac{|\nabla \phi|^2}{\phi^2}\,dx-\int_G Vu^2\,dx\,(\ge)=\int_G u^2\,\frac{f}{\phi}\,dx.
\end{eqnarray*}
This proves \eqref{A-ground} on $H^1_c(G)\cap L^\infty_c(G)$ and implies, in particular,
that $\E_V(u,u)\ge 0$ on $H^1_c(G)\cap L^\infty_c(G)$.
Therefore the extended Dirichlet space $\D(\E_V,G)$ is well defined.

Let $u\in\D(\E_V,G)$ and let
$(\varphi_n)\subset H^1_c(G)\cap L^\infty_c(G)$ be an approximating sequence for $u$.
Then $u_n:=-u\vee\varphi_n\wedge u\in H^1_c(G)\cap L^\infty_c(G)$
is also an approximating sequence for $u$ (cf. proof of Lemma \ref{l:WMP}), and
$$0\le\int_G|\nabla\left(\frac{u_n}{\phi}\right)|^2\phi^2 dx\,(\le)=
\E_V(u_n,u_n)-\int_G u_n^2\,\frac{f}{\phi}\,dx\to
\E_V(u,u)-\int_G u^2\,\frac{f}{\phi}\,dx.$$
Hence the assertion follows by Fatou's Lemma and, in the case of equality,
by standard continuity arguments.
\end{proof}

The following straightforward corollary of Lemma \ref{l:A-ground}
is crucial in our analysis of nonexistence of positive solutions
to semilinear equation \eqref{e:main'}.

\begin{corollary}\label{APA}
Assume there exists $u\in H^1_c(G)\cap L^\infty_c(G)$ such that $\E_V(u,u)< 0$.
Then equation \eqref{e:ground} has no positive super-solution.
\end{corollary}

Another interesting application of the ground state transformation is Barta's inequality.

\begin{corollary}\label{c:Barta's}
\textsc{(Barta's inequality)}
Assume that $\E_V$ is positive definite.
Then for every  $0<\phi\in H^1_{loc}(G)$ such that $(-\Delta-V)\phi\in L^1_{loc}(G)$
the following inequality holds
\begin{equation}\label{e:34}
\inf_{x\in G}\frac{(-\Delta-V)\phi}{\phi}\,\le\,
\inf_{0\neq u\in C^\infty_c(G)}\frac{\E_V(u,u)}{\|u\|_{L^2}^2}.
\end{equation}
\end{corollary}
\begin{proof}
Set $f:=(-\Delta-V)\phi$.
We may assume $f\ge 0$ in $G$ (otherwise \eqref{e:34} is trivial).
Then Lemma \ref{l:A-ground} implies that
$$\int_G |\nabla u|^2 dx-\int_G Vu^2 dx\geq \int_G u^2\frac{f}{\phi}\,dx
\geq\inf_{x\in G}\frac{f}{\phi}\int_G u^2 dx,
\qquad\forall u\in H^1_c(G)\cap L^\infty_c(G).$$
So, the assertion follows.
\end{proof}

\begin{remark}
Note that if $-\Delta-V$ admits a principal Dirichlet eigenfunction $\phi_1>0$ in $G$,
then the equality in \eqref{e:34} is attained with $\phi=\phi_1$.
\end{remark}

\begin{small}
\section*{Acknowledgments}
The authors are grateful to Vladimir Kondratiev and Zeev Sobol
for interesting and stimulating discussions.
\end{small}

\newpage

\begin{small}

\end{small}


\begin{thebibliography}{99}

\bibitem{Adimurthi}
{\sc Adimurthi, N. Chaudhuri\ and\ M. Ramaswamy},
{\it An improved Hardy-Sobolev inequality and its application},
Proc. Amer. Math. Soc. {\bf 130} (2002), 489--505.

\bibitem{Agmon}
{\sc S. Agmon},
{\it On positivity and decay of solutions of
second order elliptic equations on Riemannian manifolds},
in {\it Methods of functional analysis and theory of elliptic equations (Naples, 1982)},
19--52, Liguori, Naples, 1983.

\bibitem{Agmon-2}
{\sc S. Agmon},
Lectures on exponential decay of solutions of second-order elliptic equations:
bounds on eigenfunctions of $N$--body Schr\"odinger operators.
Princeton Univ. Press, Princeton, NJ, 1982.
%MR 85f:35019

\bibitem{Allegretto}
{\sc W. Allegretto},
{\it Criticality and the $\lambda$-property for the elliptic equations},
J. Differential Equations {\bf 69} (1987), 39--45.
%MR 88k:35010

%\bibitem{ABC}
%{\sc A.\,Ambrosetti, H.\,Brezis and G.\,Cerami},
%{\it  Combined effect of concave and convex nonlinearities in some elliptic problems}.
%{J. Funct. Anal.} {\bf 122} (1994), 519--543.

\bibitem{Bandle}
{\sc C.\,Bandle and M.\,Ess\'en},
{\it On positive solutions of Emden equations in cone-like domains},
{Arch. Rational Mech. Anal.} {\bf 112} (1990), 319--338.

\bibitem{Bandle-Levine}
{\sc C.\,Bandle and H.\,A.\,Levine},
{\it On the existence and nonexistence of global solutions of reaction--diffusion equations
in sectorial domains},
{Trans. Amer. Math. Soc.} {\bf 316} (1989), 595--622.

\bibitem{BCN}
{\sc H.\,Berestycki, I.\,Capuzzo--Dolcetta and L.\,Nirenberg},
{\it Superlinear indefinite elliptic problems and nonlinear Liouville theorems},
{Topol. Methods Nonlinear Anal.} {\bf 4} (1994), 59--78.

\bibitem{B-Veron}
{\sc M.\,-F.\,Bidaut-V\'eron}, {\it Local and global behavior of
solutions of quasilinear equations of Emden-Fowler type}, Arch.
Rational Mech. Anal. {\bf 107} (1989), 293--324.

\bibitem{B-VP}
{\sc M.-F.\;Bidaut-Veron and S.\;Pohozaev},
{\it Nonexistence results and estimates for some nonlinear elliptic problems},
J. Anal. Math. {\bf 84} (2001), 1--49.

\bibitem{Birindelli}
{\sc I.\,Birindelli and E.\,Mitidieri},
{\it  Liouville theorems for elliptic inequalities and applications},
{Proc. Roy. Soc. Edinburgh Sect. A} {\bf 128} (1998), 1217--1247.

\bibitem{Brezis-Cabre}
{\sc H.\,Brezis and X.\,Cabr\'e},
{\it Some simple nonlinear PDE's without solutions},
{Boll. Unione Mat. Ital. Sez. B} (8) {\bf 1} (1998), 223--262.

\bibitem{Brezis-Kamin}
{\sc H.\,Brezis and S.\,Kamin},
{\it Sublinear elliptic equations in $\R^N$},
{Maunscripta Math.} {\bf 74} (1992), 87-106.

%\bibitem{Caldiroli}
%{\sc P. Caldiroli\ and\ R. Musina},
%{\it On a class of two-dimensional singular elliptic problems},
%Proc. Roy. Soc. Edinburgh Sect. A {\bf 131} (2001), 479--497.

\bibitem{Chavel}
{\sc I. Chavel},
Eigenvalues in Riemannian geometry. Academic Press, Orlando, FL, 1984.

\bibitem{Davies}
{\sc E.\,B.\,Davies},
Heat Kernels and Spectral Theory.
Cambridge University Press, Cambrige-New York-Melbourn, 1989.

\bibitem{Wang}
{\sc F. Catrina\ and\ Z.--Q. Wang},
{\it On the Caffarelli-Kohn-Nirenberg inequalities:
sharp constants, existence (and nonexistence), and symmetry of extremal functions},
{Comm. Pure Appl. Math.} {\bf 54} (2001), 229--258.

\bibitem{Dupaigne}
{\sc L.\,Dupaigne},
{\it Semilinear elliptic PDE's with a singular potential},
{J. Anal. Math.} {\bf 86} (2002), 359--398.

\bibitem{Barbatis}
{\sc  S.\,Filippas, and A.\,Tertikas},
{\it Optimizing improved Hardy inequalities},
{J. Funct. Anal.} {\bf 192} (2002), 186--233.

\bibitem{FOT}
{\sc M. Fukushima, Y. \=Oshima\ and\ M. Takeda},
Dirichlet forms and symmetric Markov processes.
de Gruyter, Berlin, 1994.

\bibitem{Gidas-Spruck}
{\sc B.\,Gidas and J.\,Spruck}, {\it Global and local behavior of
positive solutions of nonlinear elliptic equations},
{Comm. Pure Appl. Math.} {\bf 34} (1981), 525--598.

\bibitem{Gilbarg}
{\sc D.\,Gilbarg and N.\,S.\,Trudinger},
Elliptic Partial Differential Equations of Second Order.
Grundlehren der Mathematischen Wissenschften 224, Springer-Verlag,
Berlin-Heidelberg-New York, 1977.


\bibitem{Hernandez}
{\sc J. Hern\'andez, F. J. Mancebo\ and\ J. M. Vega},
{\it On the linearization of some singular, nonlinear elliptic problems and applications},
Ann. Inst. H. Poincar\'e Anal. Non Lin\'eaire {\bf 19} (2002), 777--813.

%\bibitem{Kondratiev-Egorov}
%{\sc V. Kondratiev and Yu. Egorov},
%{\it On spectral theory of elliptic operators}, Birkhauser, 1996.

\bibitem{KLS}
{\sc V.\,Kondratiev, V.\,Liskevich and Z.\,Sobol},
{\it Second--order semilinear elliptic inequalities in exterior domains},
{J. Differential Equations} {\bf 187} (2003), 429--455.

\bibitem{KLS1}
{\sc V.\,Kondratiev, V.\,Liskevich and Z.\,Sobol}, {\it Positive
super-solutions to semi-linear second--order non-divergence type
elliptic equations in exterior domains}, Preprint 2004.

\bibitem{KLSU}
{\sc V.\,Kondratiev, V.\,Liskevich, Z.\,Sobol and A.\,Us},
{\it Estimates of heat kernels for a class of second--order elliptic
operators with applications to semilinear inequalities in exterior domains},
{J. London Math. Soc.} (2) {\bf 69} (2003), 107--127.

\bibitem{KLM}
{\sc V.\,Kondratiev, V.\,Liskevich and V.\,Moroz},
{\it Positive solutions to superlinear second--order divergence type elliptic
equations in cone--like domains},
{Ann. Inst. H. Poincare Anal. Non Lineaire}, to appear.

\bibitem{KLMS}
{\sc V.\,Kondratiev, V.\,Liskevich, V.\,Moroz and Z.\,Sobol},
{\it A critical phenomenon for sublinear elliptic equations in cone--like domains},
{Bull. London Math. Soc.}, to appear.

%\bibitem{LLM}
%{ V.\ Liskevich, S.\, Lyakhova and V.\,Moroz},
%{\it Positive solutions to semilinear elliptic equations
%with critical lower order terms on cone--like domains},
%in preparation.

\bibitem{Lazer}
{\sc A. C. Lazer\ and\ P. J. McKenna},
{\it On a singular nonlinear elliptic boundary-value problem},
Proc. Amer. Math. Soc. {\bf 111} (1991), 721--730.

\bibitem{Levine}
{\sc H.\,A.\,Levine},
{\it The role of critical exponents in blowup theorems},
SIAM Rev. {\bf 32} (1990), 262--288.

\bibitem{Pohozhaev}
{\sc E.\,Mitidieri and S.\,I.\,Poho\v zaev}, {\it  A priori
estimates and the absence of solutions of nonlinear partial
differential equations and inequalities (Russian)}, Tr. Mat. Inst.
Steklova {\bf 234} (2001), 1--384.

\bibitem{Murray}
{\sc J. D. Murray},
Mathematical Biology. Springer, Berlin, 1993.

\bibitem{Nachman}
{\sc A. Nachman\ and\ A. Callegari},
{\it A nonlinear singular boundary value problem in the theory of pseudoplastic fluids},
SIAM J. Appl. Math. {\bf 38} (1980), 275--281.

%\bibitem{MaR}
%{\sc Z. M. Ma\ and\ M. R\"ockner},
%{\it Introduction to the theory of (nonsymmetric) Dirichlet forms},
%Springer, Berlin, 1992.

%\bibitem{MalyZiemer}
%{\sc J. Mal\'y\ and\ W. P. Ziemer},
%{\it Fine regularity of solutions of elliptic partial differential equations},
%Amer. Math. Soc., Providence, RI, 1997.

\bibitem{Pinchover}
{\sc Y.\,Pinchover},
{\it On positive Liouville theorems and
asymptotic behavior of solutions of Fuchsian type elliptic operators},
Ann. Inst. H. Poincar\'e Anal. Non Lin\'eaire {\bf 11} (1994), 313--341.

\bibitem{Pinsky}
{\sc R.\,G.\,Pinsky},
{Positive Harmonic Functions and Diffusion}.
Cambridge Univ. Press, 1995.

\bibitem{Pohozaev-Tesei}
{\sc S. I. Pohozaev\ and\ A. Tesei},
{\it  Nonexistence of local solutions to semilinear partial differential inequalities},
Ann. Inst. H. Poincar\'e Anal. Non Lin\'eaire {\bf 21} (2004), 487--502.

\bibitem{Galaktionov}
{\sc A. A. Samarskii, V. A. Galaktionov, S. P. Kurdyumov, A. P. Mikhailov},
Blow-up in quasilinear parabolic equations. de Gruyter, Berlin, 1995.

\bibitem{Smets}
{\sc D.\,Smets and A.Tesei},
{\it On a class of singular elliptic problems with first order terms},
{Adv. Differential Equations} {\bf 8} (2003), 257--278.

%\bibitem{Ser}
%{\sc J.\,Serrin},
%{\it Entire solutions of nonlinear Poisson equations},
%Proc. London. Math. Soc. (3) {\bf 24} (1972), 348--366.

\bibitem{Serrin-Zou}
{\sc J.\,Serrin and H.\,Zou},
{\it Cauchy--Liouville and universal boundedness theorems
for quasilinear elliptic equations and inequalities},
Acta Math. {\bf 189} (2002), 79--142.

\bibitem{Terracini}
{\sc S.\,Terracini},
{\it On positive entire solutions to a class of equations with a singular coefficient and critical exponent},
Adv. Differential Equations {\bf 1} (1996), 241--264.

\bibitem{VZ}
{\sc J.\,L.\,Vazquez and E.\,Zuazua},
{\it The Hardy inequality and the asymptotic behavior of the heat equation with an inverse
square potential},
{J. Funct. Anal.} {\bf 173} (2000), 103--153.


\bibitem{Zhang-1}
{\sc Qi S.\,Zhang}, {\it  An optimal parabolic estimate and its
applications in prescribing scalar curvature on some open
manifolds with ${\rm Ricci}\geq 0$}, Math. Ann. {\bf 316} (2000),
703--731.

\bibitem{Zhang-2}
{\sc Qi S.\,Zhang}, {\it A Liouville type theorem for some
critical semilinear elliptic equations on noncompact manifolds},
Indiana Univ. Math. J. {\bf 50} (2001), 1915--1936.


\end{thebibliography}
\end{document}